\baselineskip=18pt
\def\qed{{\vrule height7pt width7pt
depth0pt}\par\bigskip}
\def\n{\noindent}
\def\Z{\ent}
\def\T{\tore}
\def\p{{\bf Proof:\ }}
\def\C{\comp}
\def\ra{\rightarrow}
\def\p.{{\bf Proof:\ }}
 \catcode`\@=11
\magnification 1200
\hsize=140mm \vsize=200mm
\hoffset=-4mm \voffset=-1mm
\pretolerance=500 \tolerance=1000 \brokenpenalty=5000

\catcode`\;=\active
\def;{\relax\ifhmode\ifdim\lastskip>\z@
\unskip\fi\kern.2em\fi\string;}

\catcode`\:=\active
\def:{\relax\ifhmode\ifdim\lastskip>\z@\unskip\fi
\penalty\@M\ \fi\string:}

\catcode`\!=\active
\def!{\relax\ifhmode\ifdim\lastskip>\z@
\unskip\fi\kern.2em\fi\string!}\catcode`\?=\active
\def?{\relax\ifhmode\ifdim\lastskip>\z@
\unskip\fi\kern.2em\fi\string?}

\frenchspacing

\newif\ifpagetitre        \pagetitretrue
\newtoks\hautpagetitre    \hautpagetitre={\hfil}
\newtoks\baspagetitre     \baspagetitre={\hfil}
\newtoks\auteurcourant    \auteurcourant={\hfil}
\newtoks\titrecourant     \titrecourant={\hfil}
\newtoks\hautpagegauche   \newtoks\hautpagedroite
\hautpagegauche={\hfil\tensl\the\auteurcourant\hfil}
\hautpagedroite={\hfil\tensl\the\titrecourant\hfil}

\newtoks\baspagegauche
\baspagegauche={\hfil\tenrm\folio\hfil}
\newtoks\baspagedroite
\baspagedroite={\hfil\tenrm\folio\hfil}

\headline={\ifpagetitre\the\hautpagetitre
\else\ifodd\pageno\the\hautpagedroite
\else\the\hautpagegauche\fi\fi}

\footline={\ifpagetitre\the\baspagetitre
\global\pagetitrefalse
\else\ifodd\pageno\the\baspagedroite
\else\the\baspagegauche\fi\fi}

\font\twbf=cmbx12\font\sc=cmcsc10

\def\date{le\ {\the\day}\ 
\ifcase\month\or Janvier\or \F\'evrier\or Mars\or Avril
\or Mai\or Juin\or Juillet\or Ao\^ut\or Septembre
\or Octobre\or Novembre\or D\'ecembre\fi\ {\the\year}}

\def\cf{{\it cf.\/}\ }    \def\ie{{\it i.e.\/}\ }
 \def\up#1{\raise 1ex\hbox{\sevenrm#1}}
\def\cqfd{\unskip\kern 6pt\penalty 500
\raise -2pt\hbox{\vrule\vbox to 10pt{\hrule width 4pt\vfill
\hrule}\vrule}\par}
\catcode`\@=12
\def \bg {\bigskip \goodbreak}
\def \sn {\nobreak \smallskip}

\def\ref#1&#2&#3&#4&#5\par{\par{\leftskip = 5em {\noindent
\kern-5em\vbox{\hrule height0pt depth0pt width
5em\hbox{\bf[\kern2pt#1\unskip\kern2pt]\enspace}}\kern0pt}
{\sc\ignorespaces#2\unskip},\
{\rm\ignorespaces#3\unskip}\
{\sl\ignorespaces#4\unskip\/}\
{\rm\ignorespaces#5\unskip}\par}}

\def\exo#1{\goodbreak\vskip 12pt plus 20pt minus 2pt
\line{\noindent\hss\bf
\uppercase\expandafter{\romannumeral#1}\hss}\nobreak\vskip
12pt }
\def \titre#1\par{\null\vskip
1cm\line{\hss\vbox{\twbf\halign
{\hfil##\hfil\crcr#1\crcr}}\hss}\vskip 1cm}

\def \eps {\varepsilon}

\def \ph{\varphi}
\def \frac#1#2{{{#1}\over {#2} }}

\def \reel{ {\rm I}\!{\rm R}}

\def \comp{ \;{}^{ {}_\vert }\!\!\!{\rm C}   }

\def \nat{ { {\rm I}\!{\rm N}} }

\def \rat{ {\rm Q}\kern-.65em {}^{{}_/ }} 

\def\ent{{{\rm Z}\mkern-5.5mu{\rm Z}}}

\def\tore{{\rm T\mkern-5.2mu\vrule height6.5 pt depth0pt
\mkern4mu}}

\def\N#1{\left\Vert #1\right\Vert }

\def\dess#1by#2(#3){\vbox to #2{
\hrule width #1 height 0pt depth 0pt
\vfill\special{picture #3}
}}

\def \proof {\nobreak \smallskip \noindent {\bf Proof.\ }}
\hfuzz=0,9cm
\centerline{\bf MULTIPLIERS AND LACUNARY SETS IN NON-AMENABLE GROUPS }
 \centerline{{by Gilles Pisier}\footnote*{Supported in
part by N.S.F. grant DMS 9003550}}
\bigskip
\bigskip
{\bf \S~0. Introduction.}

Let $G$ be a discrete group.

Let $\lambda : G \to B(\ell_2(G),\ell_2(G))$ be the left regular representation.
A function $\ph : G \to \comp$ is called a completely bounded  multiplier (= Herz-Schur
multiplier) if the transformation defined on the linear span $K(G)$ of
$\{\lambda(x),x \in G\}$ by 
$$\sum_{x \in G} f(x) \lambda(x) \to \sum_{x \in G} f(x) \ph(x) \lambda(x)$$
is completely bounded (in short c.b.) on the $C^*$-algebra $C_\lambda^*(G)$ which
is generated by $\lambda$ ($C_\lambda^*(G)$ is the closure of $K(G)$ in
$B(\ell_2(G),\ell_2(G))$.)

One of our main results  (stated below as Theorem 0.1) gives a simple
characterization of the functions $\ph$ such that $\eps \ph$ is a c.b.
multiplier on $C_\lambda^*(G)$ for any bounded function $\eps$, or equivalently
for any choice of signs $\eps(x) = \pm 1$.
We wish to consider also the case when this holds for ``almost all" choices of
signs. To make this precise, equip $\{-1,1\}^G$ with the usual uniform
probability measure. We will say that $\eps \ph$ is a c.b. multiplier of
$C_\lambda^*(G)$ for almost all choice of signs $\eps$ if   there is a
measurable subset $\Omega \subset \{-1,1\}^G$ of full measure (note that
$\Omega$ depends only on countably many coordinates) such that for any $\eps$ in
$\Omega~~~\eps \ph$ is a c.b. multiplier of $C_\lambda^*(G)$. (Note that $\ph$
is necessarily countably supported when this holds, so
the measurability issues are irrelevant.)

\proclaim Theorem 0.1. The following properties of a function $\ph : G \to
\comp$ are equivalent \sn
(i) For all bounded functions $\eps : G \to \comp$ the pointwise product $\eps
\ph$ is a c.b. multiplier.\sn
(ii) For almost all choices of signs $\eps \in \{-1,1\}^G$, the product $\eps
\ph$ is a c.b. multiplier.\sn
(iii) There is a constant $C$ and a partition of $G \times  G$ say $G \times  G
= \Gamma_1 \cup \Gamma_2$ such that 
$$\sup_{s \in G} \sum_{t \in G} |\ph(st)|^2 1_{\{(s,t) \in \Gamma_1\}} \leq
C^2~{\rm and}~\sup_{t \in G} \sum_{s \in G} |\ph(st)|^2
1_{\{(s,t) \in \Gamma_2\}} \leq C^2.$$\sn
(iv) There is a constant $C$ such that for all finite subsets $E,F \subset G$
with $|E| = |F|= N$ we have
$$\sum_{(s,t) \in E \times  F} |\ph (st)|^2 \leq C^2 N.$$\sn
(v) There is a constant $C$ such that for any Hilbert space $H$ and for any
finitely supported function $a : G \to B(H)$ we have
$$\N {\sum_{x \in G} \ph(x)\lambda(x) \otimes
a(x)}_{B(\ell_2(G,H))} \leq C \max \left\{\N {\left(\sum
a(x)^* a(x)\right)^{1/2}},\N {\left(\sum a(x)
a(x)^*\right)^{1/2}}\right\}.$$

 {\bf Note: } The properties (iii) and (iv) could have been
stated equivalently with the function
$(s,t) \to \ph(st^{-1})~$ (which would have been
perhaps more natural)${\rm or}~(s,t) \to \ph(s^{-1}
t)$ instead of $(s,t) \to \ph(st)$. We chose the simplest
notation.

This theorem is proved in section 2 below.

\n{\bf Remark: } The papers [W] and [B3] show that amenable groups are
characterized by the property that all multipliers $\ph$ satisfying
(iii) in Theorem 0.1 are necessarily in $\ell_2(G)$. Hence
the preceding statement is of interest only in the non-amenable case.
Moreover, on the free
group with finitely many generators, the radial functions 
which satisfy the above property (iii) are characterized in
[W].

The equivalence of (iii) and (iv) is already known. It was proved by Varopoulos
[V1] in his study of the projective tensor product
$\ell_\infty \hat \otimes \ell_\infty$ and the Schur
multipliers of $B(\ell_2,\ell_2)$. Let $(e_s)$ (resp.
$(e_t)$) denote the canonical basis of $\ell_\infty(S)$
(resp. $\ell_\infty(T)$). We will denote by $\tilde V(S,T)$
the set of functions $\psi : S \times  T \to \comp$ such
that $$ \sup_{\buildrel{E \subset S,F \subset T}\over {|E|<
\infty,|F|< \infty}} \N {\sum_{\buildrel{S \in E}\over {t
\in F}} \psi (s,t) e_s \otimes e_t}_{{\ell_\infty} (S) \hat
\otimes {{\ell_\infty}(T)}} < \infty.$$

More precisely, Varopoulos proved

\proclaim Theorem 0.2. ([V1]) Let $S,T$ be arbitrary sets. The following
properties of a function $\psi : S \times  T \to \comp$ are equivalent.\sn
(i) For all bounded functions $\eps : S \times  T \to \comp$ the pointwise
product $\eps \psi$ is in $\tilde V(S,T)$.\sn
(ii) For almost all choices of signs $\eps$ in $\{-1,1\}^{S\times T}$ the
pointwise product $\eps \psi$ is in $\tilde V(S,T)$.\sn
 (iii) There is a constant $C$ and a partition $S
\times  T = \Gamma_1 \cup \Gamma_2$ such that
$$\sup_{s \in S} \sum_{\buildrel{t \in T}\over {(s,t) \in \Gamma_1}} |\psi
(s,t)|^2 \leq C^2~{\rm and}~\sup_{t \in T}
\sum_{\buildrel{s \in S}\over {(s,t) \in \Gamma_2}} |\psi
(s,t)|^2 \leq C^2.$$\sn (iii)' There is a decomposition
$\psi = \psi_1 + \psi_2$ with  $$\sup_{s \in S} \sum_{t \in
T} |\psi_1 (s,t)|^2 < \infty~{\rm and}~\sup_{t \in T}
\sum_{s \in S} |\psi_2(s,t)|^2 \leq C.$$\sn
(iv) There is a constant $C$ such that for all finite
subsets $E \subset S,F \subset T$ with $|E| = |F| = N$, we
have  $$\sum_{(s,t) \in E \times F} |\psi(s,t)|^2 \leq C^2
N.$$\sn

The deepest implication in Theorem 0.2 is (ii) $\Rightarrow$
(iii). The equivalence of (iii) and (iii)' is obvious and
(iii)' $\Rightarrow$ (i) is rather easy (by duality, it
follows from Khintchine's inequality). The equivalence (iii)
$\Leftrightarrow$ (iv) is a remarkable fact of independent
interest. The decompositions of the
form (iii) are related to some early work of Littlewood and
the matrices admitting the decomposition (iii) are often
called Littlewood tensors, following Varopoulos's
terminology. We note in passing that (ii) $\Rightarrow$
(iii) (and in fact a slightly stronger result) can be
obtained as an application of Slepian's comparison
principle for Gaussian processes in the style of S. Chevet
(see [C] th\'eor\`eme 3.2). However, we do not see how to
exploit this approach in our more general context.

We will prove below a result which contains Theorem 0.2 as a particular case
and implies Theorem 0.1 in the group case. Roughly our result gives a necessary
condition (analogous to the above (iii)) for a random
series $\sum_{n=1}^\infty \eps_n \psi_n$ with random signs
$\eps_n = \pm 1$ and arbitrary coefficients $\psi_n$ in
$\ell_\infty \hat \otimes \ell_\infty$ to define a.s. an
element of  $\ell_\infty \hat \otimes \ell_\infty$. The
implication (ii) $\Rightarrow$ (iii) in Theorem 0.2
corresponds to the particular case when $\psi_n$ is of the
form $\psi_n = \alpha_n\  e_{i_n} \otimes e_{j_n}$ where
$\alpha_n \in \comp$ and $n \to (i_n,j_n)$ is a bijection
of $\nat$ onto $\nat \times  \nat$.

Our necessary condition can be stated as follows : there is a sequence of
scalars $\alpha_m$ with $\sum_{m \geq 0} |\alpha_m|< \infty$ and scalar
coefficients
$$a_n^m(i),b^m(j),c^m(i),d_n^m(j)$$
such that
$$\psi_n(i,j) = \sum_{m \geq 0} \alpha_m [a_n^m(i) b^m(j) + c^m(i) d_n^m(j)]$$
and such that for all $m$
$$\sup_i (\sum_n |a_n^m(i)|^2)^{1/2} \leq 1,\quad \sup_j|b^m(j)| \leq 1$$
$$\sup_i |c^m(i)|\leq 1,\quad \sup_j (\sum_n |d_n^m (j)|^2)^{1/2} \leq 1.$$
In other words, the condition expresses that the sequence $(\psi_n)$ can be
written (up to a multiplicative norming constant) as an element of the
closed convex hull of special sequences of the form
$$\psi_n(i,j) = a_n(i) b(j) + c(i) d_n(j)$$
with
$$\sum_n |a_n(i)|^2 \leq 1,|b(j)| \leq 1, |c(i)|\leq 1, \sum_n |d_n(j)|^2 \leq
1$$
for all $i$ and $j$.

This will be stated below (\cf Theorems 2.1 and 2.2) in the more precise (and
concise) language of tensor products. 

To emphasize the content of Theorem 0.1, we now state an application in terms
of Schur multipliers. For any sets $S,T$ a function $\psi : S \times  T \to
\comp$ is called a Schur multiplier of $B(\ell_2(S),\ell_2(T))$ if for any $u
\in B(\ell_2(S),\ell_2(T))$ with associated matrix $(u(s,t))$ the matrix
$(\psi(s,t) u(s,t))$ is the matrix of an element of $B(\ell_2(S),\ell_2(T))$.
It is known that the set of all Schur multipliers $\psi : S \times  T \to
\comp$ coincides with the space $\tilde V(S,T)$. This essentially goes back to
Grothendieck [G]. We give more background on Schur multipliers in section 1. The
next statement is an application of Theorem 0.1 (and the easier implication
(iii) $\Rightarrow$ (i) in Theorem 0.2).

\proclaim Corollary 0.3. Assume that $\ph$ satisfies (i) in
Theorem 0.1. Then for all choices of signs $\xi \in
\{-1,1\}^{G \times  G}$ (indexed by $G \times G$ this time)
the product $$(s,t) \to \xi (s,t) \ph(st)$$
is in $\tilde V(G,G)$ hence it defines a Schur
multiplier of $B(\ell_2(G),\ell_2(G))$.

Actually, the group structure plays a rather limited role in the preceding
statement and in Theorem 0.1. To emphasize this point we
state (see also Remark 2.4 below)

\proclaim Corollary 0.4. Let $G$ be any set.
Suppose given a map
$$p : G \times  G \to G$$
such that for all fixed $(s_0,t_0)$ in $G \times  G$ the maps $s \to p(s,t_0)$
and $t \to p(s_0,t)$ are bijective. (Actually it suffices to assume that there
is a fixed finite upper bound on the cardinality of the sets $\{s |p(s,t_0) =
x\}$ and $\{t |p(s_0,t) = x\}$ when $x,s_0,t_0$ run over
$G$). Let $\ph : G \times  G \to \comp$ be a function on $G
\times  G$. Assume that for all (actually ``almost all" is
enough) choices of signs $(\eps_x)_{x \in G}$ the function
$$(s,t) \to \eps_{p(s,t)} \ph(p(s,t))$$
is a Schur multiplier of $B(\ell_2(G),\ell_2(G))$. Then, for all choices of
signs $\eps_{s,t}$ (indexed by $G \times  G$ this time) the function
$$(s,t) \to \eps_{s,t} \ph (p(s,t))$$
is a Schur multiplier of $B(\ell_2(G),\ell_2(G))$.

The results stated above are proved in section 2. In section 3, we apply them
to study a class of ``lacunary subsets" of a discrete group which is analogous
of the class of finite unions of Hadamard-lacunary subsets of $\nat$. We give a
combinatorial characterization of these sets which we call $L$-sets, but we
leave as a conjecture a stronger result (see conjecture 3.5 below).

\vfill\eject

{\bf \S~ 1. Preliminary Background.}

We refer to [Pa] for more information on completely bounded maps.

Let $S$ be any set. As usual we denote by $\ell_\infty(S)$ the space of all
complex valued bounded functions on $S$, equipped with the sup-norm.

For any Banach space $E$, we will also use the space $\ell_\infty(S,E)$ of all
$E$-valued bounded functions $x : S \to E$ equipped with the norm $\N {x} =
\sup_{s \in S} \N {x(s)}_E$.

When $S = \nat$, we write simply $\ell_\infty$. In particular, we will use below
the space $\ell_\infty (\ell_2)$ which also can be
regarded  as the space of all matrices $x(j,k)$ such that
$$\sup_j (\sum_k |x (j,k)|^2)^{1/2}< \infty.$$

Let $X,Y$ be Banach spaces and let $X \otimes Y$ be their linear tensor product.
We   recall the definition of the projective norm and of several other
important tensor norms (\cf [G]).

For any $u$ in $X \otimes Y$, let 
$$\N {u}_\wedge = \inf \{\sum_1^n \N {x_i} \N {y_i}| u = \sum_1^n x_i \otimes
y_i,\ x_i \in X,\  y_i \in Y\}\leqno (1.1)$$

We will also need
$$\gamma_2(u) = \inf\{\sup_{\xi \in B_{X^*}} (\sum_1^n |\xi(x_i)|^2)^{1/2}
\sup_{\eta \in B_{Y^*}}(\sum_1^n |\eta (y_i)|^2)^{1/2}\}\leqno (1.2)$$
where the infimum runs again over all possible representations of the form $u =
\sum_1^n x_i \otimes y_i$.

Equivalently $\gamma_2(u)$ is the ``norm of factorization through a Hilbert
space" of the associated operator $u : X \to Y^*$.

We will also need a generalization of the $\gamma_2$-norm considered in [K] for
Banach lattices. Recall that a Banach lattice $X$ is called $2$-convex if we
have
$$\forall~x,y \in X~~~~~\N {(|x|^2 + |y|^2)^{1/2}} \leq (\N {x}^2 + \N
{y}^2)^{1/2},$$
see e.g. [LT] for more information.

Let $X,Y$ be two $2$-convex Banach lattices. For $u = \sum_1^n x_i \otimes y_i
\in X \otimes Y$, we define 
$$\gamma(u) = \inf \{\N {(\sum_1^n |x_i|^2)^{1/2}}_X \N {(\sum_1^n
|y_i|^2)^{1/2}}_Y\}\leqno (1.3)$$
where the infimum runs over all representations of $u$.

It is easy to see that the $2$-convexity of $X$ and $Y$ implies that this is a
norm on $X \otimes Y$. Note that $\ell_\infty$ (or more
generally $L_p$ for $2\leq p\leq \infty)$ is an example of a
$2$-convex Banach lattice. In the case of the product
$\ell_\infty \otimes \ell_\infty$ it is easy to check that
(1.2) and (1.3) are identical, so that  $$\gamma = \gamma_2
~{\rm on}~\ell_\infty \otimes \ell_\infty \leqno (1.4)$$

Indeed, if $x_1,...,x_n \in \ell_\infty(S)$ over an index set $S$ we have
clearly
$$\sup_{\xi \in B_{\ell_\infty^*}} (\sum |\xi(x_i)|^2)^{1/2} =
\sup_{\sum|\alpha_i|^2 \leq 1} \N {\sum \alpha_i
x_i}_{\ell_\infty} = \sup_{s \in S} (\sum
|x_i(s)|^2)^{1/2}=\{\N {(\sum 
|x_i|^2)^{1/2}}_{\ell_\infty}$$

Let $S$ and $T$ be two index sets. Consider $u$ in $\ell_\infty (S) \otimes
\ell_\infty(T)$ with associated matrix $u(s,t) = < \delta_s \otimes \delta_t,u >$
(we denote by $(\delta_s)$ and $(\delta_t)$ the Dirac masses at $s$ and $t$
respectively, viewed as linear functionals on 
$\ell_\infty(S)$ and $\ell_\infty(T))$. Then we have $\gamma_2(u) \leq 1$ iff
there are maps $x : S \to \ell_2~~~~~y : T \to \ell_2~~~{\rm such ~that}$
$\sup_{s \in S} \N {x(s)} \leq 1,\  \sup_{t \in T} \N {y(t)} \leq 1~{\rm and}$
$$\forall ~s,t \in T~~~~u(s,t) = < x(s),y(t) >.$$

This is very easy to check.

The following result is well known 

\proclaim Proposition 1.1. Let $S,T$ be arbitrary sets. Let $\ph : S \times T
\to \comp$ be a function.  We consider the Schur multiplier
$$M_\ph : B(\ell_2(S),\ell_2(T)) \to B(\ell_2(S),\ell_2(T))$$
defined in matrix notation by $M_\ph ((a(s,t))) = (\ph(s,t) a(s,t))$.
The following are equivalent\sn
(i) $\N {M_\ph} \leq 1$,\sn
(ii) There are vectors $x(s),y(t)$ in a Hilbert space   such
that\sn
 $\sup_{s}\N {x(s)} \leq 1, \sup_t \N {y (t)} \leq 1$ and $\ph(s,t) = <
x(s),y(t) >$.\sn
(iii) For all finite subsets $E \subset S$ and $F \subset T$ we have
$$\N {\sum_{(s,t) \in E \times  F} \ph(s,t) e_s \otimes e_t}_{\ell_\infty(S)
\hat \otimes_{\gamma_2} \ell_\infty(T)} \leq 1.$$
Moreover if $S$ and $T$ are finite sets then (i) (ii) and
(iii) are equivalent to\sn (iv)~~~~~$\N {\sum_{(s,t) \in
S\times T} \ph(s,t) e_s \otimes e_t}_{\ell_\infty(S) \hat
\otimes_{\gamma_2} \ell_\infty (T)} \leq 1$ where $(e_s)$
and $(e_t)$ denote the canonical bases of $\ell_\infty(S)$
and $\ell_\infty(T)$ respectively.

\proof Let us first assume that $S$ and $T$ are \underbar
{finite} sets. The equivalence of (ii) (iii) and (iv) is then obvious. Assume
(i). This means exactly that for any $a : \ell_2(S) \to \ell_2(T)$ with $\N {a}
\leq 1$ and for any $\alpha$ and $\beta$ in the unit ball respectively of
$\ell_2(S)$ and $\ell_2(T)$ we have
$$|\sum_{s,t} \ph(s,t) a(s,t) \alpha(s) \beta(t)|\leq 1.$$

In other words $\N {M_\ph}\leq 1$ means that $\ph$ lies in the polar of the
set $C_1$ of all matrices of the form $(\alpha(s) a(s,t) \beta(t))$ with $a,
\alpha, \beta$ as above.
But it turns out that this set $C_1$ is itself the polar of the set $C_2$ of
all matrices $(\psi (s,t))$ such that $\N {\sum \psi (s,t) e_s \otimes
e_t}_{\ell_\infty (S) \otimes_{\gamma_2}\ell_\infty (T)} \leq 1$. (Indeed, this
follows from the known factorization property which describes the norm
$\gamma^*_2$ which is dual to the norm $\gamma_2$,cf.
e.g. [Kw] or [P1] chapter 2.b). In conclusion $\ph$ belongs
to $C_2^{00} = C_2$ iff (i) holds, and this proves the
equivalence of (i) and (iv) in the case $S$ and $T$ are
finite sets.

In the general case of arbitrary sets $S$ and $T$, we note that (ii)
$\Rightarrow$ (iii) $\Rightarrow$ (i) is obvious by passing to finite subsets. It
remains to prove (i) $\Rightarrow$ (ii), but this is immediate by a compactness
argument. Indeed, if (i) holds there is obviously a net $(\ph^i)$ tending to
$\ph$ pointwise and formed with finitely supported functions on $S \times  T$
such that $\N {M_{\ph^i}} \leq 1$. Then by the first part of the proof, each
$\ph^i$ satisfies (ii) and it is easy to conclude by an ultraproduct argument
that $\ph$ also does. \qed

 \n{\bf Remark}. As observed by Uffe Haagerup (see [H3])
 Proposition 1.1 implies that the completely bounded norm of $M_\ph$
coincides with its norm. Indeed, it is easy to deduce from (ii) that $\N
{M_\ph}_{cb} \leq 1$.

In the harmonic analysis literature, the c.b. multipliers of $C_\lambda(G)$ are
sometimes called Herz-Schur multipliers. They were considered by Herz (in a
dual framework, as multipliers on $A(G)$) before the notion of complete
boundedness surfaced. The next result from [BF]
(see also [H3]) clarifies the relation between the various
kinds of multipliers.

\proclaim Proposition 1.2.   Let $G$ be a discrete group.
Consider a function $\ph : G \to \comp$. We define then
complex functions $\ph_1,\ph_2$ and $\ph_3$ on $G \times 
G$ by setting $$\forall~(s,t) \in G \times 
G~~~~\ph_1(s,t) = \ph(st^{-1}),\ph_2(s,t) = \ph(s^{-1} t),
\ph_3(s,t) = \ph(st).$$ We consider the corresponding Schur
multipliers $M_{\ph_1},M_{\ph_2}$ and $M_{\ph_3}$ on
$B(\ell_2(G),\ell_2(G))$. Then\sn (i) ([BF]) The Schur
multiplier $M_{\ph_1}$ is bounded iff the linear operator
$T_\ph : C^*_\lambda(G) \to C^*_\lambda(G)$ which maps
$\lambda(x)$ to $\ph(x) \lambda(x)$ is completely bounded.
Moreover, $\N {M_{\ph_1}} = \N {T_\ph}_{cb}$. \sn
(ii) Moreover,
$M_{\ph_1}$ is bounded iff $M_{\ph_2}$ (resp. $M_{\ph_3}$)
is bounded and we have  $$\N {M_{\ph_1}} = \N {M_{\ph_2}} =\N
{M_{\ph_2}}_{cb}= \N {M_{\ph_3}}=\N {M_{\ph_3}}_{cb}.$$

\proof The last assertion is immediate (note that if $\psi(s,t)$ is a bounded
Schur multiplier on $S\times T$ then for any bijections $f : S \to S$ and $g :
T \to T~~~~~\psi(f(s),g(t))$ also is  a bounded
Schur multiplier  with the same norm).
Note that $M_{\ph_1}$ leaves $C^*_\lambda(G)$ invariant and its restriction to
$C^*_\lambda(G)$ coincides with $T_\ph$. Hence by the preceding remark we have
$$\N {T_\ph}_{cb} \leq \N {M_{\ph_1}}_{cb} = \N {M_{\ph_1}}.$$
Conversely, if $\N {T_\ph}_{cb} \leq 1$ then the factorization theorem of c.b.
maps due to Wittstock (Haagerup [H3] and Paulsen proved it
independently, see [Pa]) says that there is a Hilbert space
$H$, a representation $\pi : B(\ell_2(G)) \to B(H)$ and
operators $V_1$ and $V_2$ from $\ell_2(G)$ into $H$ with
$\N {V_1} \leq 1,\N {V_2} \leq 1$ such that $\forall~a \in
C^*_\lambda(G)~~~~T_\ph(a) = V^*_2 \pi(a) V_1$.

In particular we have $\ph(x) \lambda(x) = T_\ph(\lambda(x)) = V^*_2
\pi(\lambda(x)) V_1$, which implies 
$$\eqalign{\forall~s,t \in G~~~~~\ph(st^{-1}) & =  <
\delta_s,T_\ph(\lambda(st^{-1})) \delta_t >\cr
& = < \pi(\lambda(s))^* V_2 \delta_s, \pi(\lambda(t^{-1})) V_1 \delta_t>}$$
This shows that $\ph_1$ satisfies (ii) in Proposition 1.1, hence $\N
{M_{\ph_1}}\leq 1$. \qed

Grothendieck [G] proved that $\gamma_2$ and $\N {~}_\wedge$ are equivalent
norms on $\ell_\infty \otimes \ell_\infty$ (or on $\ell_\infty(S) \otimes
\ell_\infty (T)$, more precisely there is a constant $K_G$ such that
$$\forall~u \in \ell_\infty \otimes \ell_\infty ~~~~~~~~~\N {u}_\wedge \leq
K_G \gamma_2(u) .\leqno (1.5)$$
The exact numerical value of the best constant $K_G$ in (1.5) is still an open
problem (see [P1] for more recent results). 

Grothendieck's striking theorem admits many equivalent reformulations. In the
context of Banach lattices, Krivine [K] emphasized the following one. Let $X,Y$
be 2-convex Banach lattices, then $\gamma$ and $\N {~}_\wedge$ are equivalent
norms on $X \otimes Y$ and we have
$$\forall ~u \in X \otimes Y~~~~~~~~~\N {u}_\wedge \leq K_G \gamma (u).
\leqno (1.6)$$
Note that in (1.5) and (1.6) the converse inequality is trivial (since $\N
{~}_\wedge$ is the ``greatest cross-norm"), we have $\gamma_2 (u) \leq \N
{u}_\wedge$ and $\gamma (u) \leq \N {u}_\wedge$ for all $u$ in $X \otimes Y$. 

The reader should note that the equality $\gamma = \gamma_2$ on $\ell_\infty
\otimes \ell_\infty$ is a special property of $\ell_\infty$ spaces. If $X = Y =
\ell_2$ for instance then on $X \otimes Y~~~~\gamma_2$ is the injective norm
(\ie the usual operator norm) while $\gamma$ is identical to the projective
norm (\ie the trace class norm).

We refer the reader to [P2] for the discussion of a more general class of
cross-norms which behave like $\gamma$ and $\gamma_2$.

While the proof of Proposition 1.1 uses nothing more than
the Hahn-Banach theorem, the next result is a reformulation
of Grothendieck's theorem one more time, it was observed in
some form already in [G] (Prop. 7, p. 68), and was later
rediscovered and extended by various authors, notably
J.Gilbert in harmonic analysis (see [GL],[Be]) and
U.Haagerup in operator algebras (see [H3] the unpublished
preliminary version of [H2]).

\proclaim Theorem 1.3. In the case when $S,T$ are finite sets in the same
situation as Proposition 1.1 we have 
$${1 \over K_G} \N {\sum \ph(s,t) e_s \otimes e_t}_{\ell_\infty (S) \hat \otimes
\ell_\infty(T)} \leq \N {M_\ph}\leq \N {\sum \ph(s,t) e_s \otimes
e_t}_{\ell_\infty (S) \hat \otimes \ell_\infty(T)}.$$
Moreover, when $S,T$ are arbitrary  sets, the space of all
bounded Schur multipliers of $B(\ell_2(S),\ell_2(T))$
coincides with the space $\tilde{V}(S,T)$.

\proof This follows immediately from Grothendieck's
inequality (1.5) and Proposition 1.1.\qed

Perhaps a more
intuitive formulation is as follows. Let us call ``simple
multipliers" the Schur multipliers of the form  $$\ph(s,t)
= \eps_s \eta_t$$ with $\eps_s,\eta_t \in \comp$ such that
$|\eps_s| \leq 1,|\eta_t| \leq 1$. These are obviously such
that $\N {M_\ph} \leq 1$, but precisely Theorem 1.3 says
that any multiplier $\ph$ with $\N {M_\ph} \leq {1 \over
K_G}$ lies in the
  convex hull of the set of simple multipliers if $S,T$ are finite sets and if
$S,T$ are infinite sets, then $\ph$ lies in the pointwise closure of the
  convex hull of the set of simple multipliers.

  \n{\bf Remark 1.4: } Consider again a ``simple Schur multiplier" of the
form $\ph(s,t) = \eps_s \eta_t$ as above with $|\eps_s|\leq 1, |\eta_t| \leq
1$. Then we have
$$\forall~A \in B(\ell_2(S),\ell_2(T))~~~~~~~M_\ph(A) = v_1 A v_2$$
where $v_1 : \ell_2(S) \to \ell_2(S)$ and $v_2 : \ell_2(T) \to \ell_2(T)$ are
the diagonal operators of multiplication by $(\eps_s)$ and $(\eta_t)$
respectively. Therefore, it is obvious that $\N {M_\ph}_{cb} \leq 1$.
(see [Pa] for more information.)

We need to consider ``sums" of Banach spaces which are usually not direct sums.
Although we will mainly work with natural concrete Banach spaces $X$ and $Y$
for which saying that an element belongs to $X + Y$ will have a clear meaning,
we recall the following formal definition of $X + Y$.

Assume that $X,Y$ are both continuously injected in a larger topological vector
space ${\cal X}$. Then $X + Y$ is defined as the subspace of ${\cal X}$ of all
elements of the form $\sigma = x + y$ with $x \in X, y \in Y$, equipped with
the norm
$$\N {\sigma}_{X + Y} = \inf \{\N {x}_X + \N {y}_Y |\sigma = x + y\}.$$
Equipped with this norm, $X + Y$ is a Banach space (and its dual can be
identified with $X^* \cap Y^*$ under some mild compatibility assumption on
$X,Y$).

Alternately, one may consider the direct sum $X \oplus_1 Y$ equipped with the
norm $\N {(x,y)} = \N {x} + \N {y}$ together with the closed subspace $N
\subset X \oplus_1 Y$ of all elements $(x,y)$ which satisfy the identity
$x+y=0$ when injected into ${\cal X}$.

Then the quotient space $\sum = (X \oplus_1 Y)/N$ can be identified with $X +
Y$.

It will be convenient at some point to use the following elementary fact.

\proclaim Lemma 1.5. Let $dm(t) = {dt \over 2 \pi}$ be the normalized
Haar measure on $\tore$. Then for any integer $N$ and any continuous function
$f : \tore^N \to \reel$ we have
$$\int f(t_1,...,t_N) dm (t_1)... dm(t_N) \geq \inf \int f(e^{in_1 t},e^{in_2
t} ,...,e^{in_N t}) dm(t)$$
where the infimum on the right side runs over all sets of integers
$n_1,n_2,...,n_N$ with $2 n_1 < n_2,2n_2 < n_3,...,2n_{N-1} < n_N$.

\proof If $f$ is a trigonometric polynomial, this is obvious by choosing $(n_k)$
lacunary enough. By density, this must remain true for all real valued $f$ in
$C(\tore^N)$.
Let us consider the infinite dimensional torus ${\T}^{\nat}$. We denote by
$z_j$ the $j$-th coordinate on ${\T}^{\nat}$ and by $\mu$ the normalized Haar measure on
${\T}^{\nat}$.
The following is a reformulation of the main result of [LPP].

\proclaim Theorem 1.6. Let $a_1,...,a_n$ be elements of a von Neumann algebra
$M$, let $\xi_1,...,\xi_n$ be elements of the predual $M_*$. Then 
$$|\sum_1^n < \xi_j,a_j>| \leq \int \N {\sum z_j a_j}_{M_*} d\mu(z) [\N {(\sum
a^*_j a_j)^{1/2}}_M + \N {(\sum a_j a_j^*)^{1/2}}_M].\leqno (1.7)$$

\proof Two approaches are given in [LPP]. The first one proves this result using
the factorization of analytic functions in $H^1$ with values in $M_*$. Actually,
in [LPP] (1.7) is stated with $\int \N {\sum z_j \xi_j}_{M_*} d\mu (z)$
replaced by $\int \N {\sum e^{in_jt} \xi_j}_{M_*} {dm(t) }$ for any
lacunary sequence $n_j$ such that $n_j > 2 n_{j-1}$. Using the preceding lemma,
it is then easy to obtain  (1.7) as stated above. (Moreover, it is possible
to use the factorization argument of [LPP] 
directly in $\tore^N$, see the following remark.)
A second approach is given in the appendix of [LPP]. There it is shown that
(1.7), with some additional numerical factor, can be deduced from (and is
essentially equivalent to) the non-commutative Grothendieck inequality due to
the author (see [P1], Theorem 9.4 and Corollary 9.5).

\n{\bf Remark}.The reader may find
 the use of a lacunary sequence $(n_k)$ in the preceding
proof a bit artificial. Actually, we can use directly the independent sequence
$(z_k)$ on $\tore^\nat$ equipped with $\mu$. Indeed, the classical
factorization theory of $H^1$ functions as products of two $H^2$ functions
extends to this setting, provided one considers $\tore^\nat$ as a compact group
with ordered dual in the sense e.g. of [R] chapter 8. Here the dual of
$\tore^\nat$  is ordered lexicographically.
The factorization of matrix valued functions (as used in [P2] Appendix B) also
extends to this setting, so that the main results of [P2] also remain valid in
this setting. This approach is described in [P3]. We chose the more traditional
``one dimensional" torus presentation to provide more precise and explicit
references for the reader.

We will use the following well known consequence of the Hahn-Banach Theorem
(\cf [Kw], see also e.g. Lemma 1.3 in [P2]).

\proclaim Lemma 1.7. Let $S,T$ be finite sets.  
Let $u : \ell_\infty(S,\ell_2) \hat \otimes_\gamma \ell_\infty (T) \to \comp$
be a linear form of norm $\leq 1$ on $\ell_\infty(S,\ell_2) \hat \otimes_\gamma
\ell_\infty (T)$. Then there are probabilities $P,Q$ on $S$ and $T$ such that
$\forall~\ph \in \ell_\infty (S,\ell_2)~~~\forall~\eta \in \ell_\infty (T)$.
$$|< u,\ph \otimes \eta>| \leq (\int \N {\ph(s)}^2_{\ell_2} dP(s))^{1/2}
(\int |\eta(t)|^2 dQ (t))^{1/2}. \leqno (1.8)$$

\proof (Sketch). By assumption and by definition (1.3) we have for all finite
sequences $\ph_k \in \ell_\infty (S,\ell_2),\eta_k \in \ell_\infty(T)$
$$\eqalign{\sum_k|< u,\ph_k \otimes \eta_k >|& \leq \sup_{s \in S} (\sum \N
{\ph_k(s)}^2)^{1/2} \sup_{t \in T} (\sum |\eta_k(t)|^2)^{1/2}\cr
& \leq {1 \over 2} \sup_{S \times T} \{\sum \N {\ph_k (s)}^2 + |\eta_k
(t)|^2\}}$$
Let $C$ be the convex cone in $C(S \times  T)$ formed by all the functions of
the form
$$(s,t) \to {1 \over 2} \sum \N {\ph_k(s)}^2 + |\eta_k(t)|^2 - |< u,\ph_k
\otimes \eta_k >|.$$
Then $C$ is disjoint from the open cone $C_{-} = \{ \ph | \max~ \ph < 0\}$,
hence (Hahn-Banach) there is a hyperplane in $C (S \times T)$ which separates
$C$ and $C_{-}$.
By an obvious adjustment, this yields a probability $\lambda$ on $S \times  T$
such that $\int f(s,t) d \lambda \geq 0$ for any $f$ in
$C$. Hence letting $P$ (resp. $Q$) be the projection of
$\lambda$ on the first (resp. second) coordinate we obtain
$$|< u,\ph \otimes \eta>|\leq {1 \over 2} (\int \N {\ph(s)}^2 dP(s) + \int
|\eta (t)|^2 dQ(t)).$$
Finally applying this to $(\ph \theta^{-1}) \otimes (\theta \eta)$ and
minimizing the right hand side over all $\theta > 0$, we obtain the announced
result (1.8). \qed

\n{\bf Remark 1.8 }
 Let $S,T$ be finite sets and let   $P,Q$ be probabilities  on $S$ and
$T$ respectively. Let $J_P : \ell_\infty(S) \to L_2(P)$ and
$J_Q : \ell_\infty(T) \to L_2(Q)$ be the canonical
inclusions.
 Then we have 
   $\forall~\psi \in \ell_\infty(S) \otimes
\ell_\infty(T)$ $$\N {(J_P \otimes J_Q)(\psi)}_{L_2(P) \hat
\otimes L_2(Q)} \leq \N {\psi}_{\ell_\infty(S) \hat
\otimes_{\gamma_2} \ell_\infty(T)}. \leqno (1.9)$$
Indeed, this is elementary. For any $x_1,...,x_n$ in
${\ell_\infty(S)}$, $y_1,...,y_n$ in ${\ell_\infty(T)}$,
we have 
$$\sum\| x_i\|^2_{L_2(P)})^{1/2} \leq
\|(\sum|x_i|^2)^{1/2}\|_{\ell_\infty(S)}\quad {\rm
and}\quad \sum\| y_i\|^2_{L_2(P)})^{1/2} \leq
\|(\sum|y_i|^2)^{1/2}\|_{\ell_\infty(T)}.$$
This clearly implies (1.9). 

\baselineskip=18pt
\def\qed{{\vrule height7pt width7pt
depth0pt}\par\bigskip}
\def\Z{\ent}
\def\T{\tore}
\def\n{\noindent}
\def\p{{\bf Proof:\ }}
\def\C{\comp}
\def\ra{\rightarrow}
 
\n {\bf Remark 1.9}  Let
 us denote simply by $H_1(\T; M_*)$ the subspace of $L_1(\T,dm;M_*)$ formed of all
the functions $f$ such that the ($M_*$-valued ) Fourier transform 
is supported on the non-negative integers. 
Similarly, we can denote by $H_1({\T}^{\nat}; M_*)$ the subspace of
$L_1({\T}^{\nat},{m}^{\nat};M_*)$ formed by
 the functions with Fourier transform supported by
the non-negative elements of ${\Z}^{(\nat)}$ ordered
  lexicographically. We again denote by $z_j$ the $j$-th coordinate on
${\T}^{\nat}$ and we let $\hat{f}(z_j) =\int f \bar{z}_j$.
In [LPP] the following refinement of (1.7) is proved.

Assume 
  that there is a   function $f $ in the unit ball of $H_1(\T; M_*)$ (resp.
$H_1({\T}^{\nat}; M_*)$) 
such that $\hat{f}(3^j)=a_j$ (resp. $\hat{f}(z_j) =a_j$) for all $j$ ,
 then we have 
$$|\sum_1^n < \xi_j,a_j>| \leq   [\N {(\sum
a^*_j a_j)^{1/2}}_M + \N {(\sum a_j a_j^*)^{1/2}}_M].\leqno (1.10)$$ 

\vfill\eject
  {\bf \S~2. Main results.}

Our main result is a general statement which does not use the group structure
at all, it can be viewed as a generalization of Varopoulos's result stated
above as Theorem 0.2.

\proclaim Theorem 2.1. Let $S$ and $T$ be arbitrary sets, and let
$$\psi_n \in \ell_\infty(S) \hat \otimes \ell_\infty(T)$$
be a sequence such that the series
$$\sum_{n=1}^\infty \eps_n \psi_n$$
converges in
$\ell_\infty(S) \hat \otimes \ell_\infty(T)$ for almost all choice of signs
$\eps_n = \pm 1$. Then, if we denote by $(e_n)$ the canonical basis of
$\ell_2$, the series 
$$\sum_{n=1}^\infty e_n \otimes \psi_n$$
is convergent in the space $\ell_\infty(S,\ell_2) \hat \otimes \ell_\infty(T) +
\ell_\infty(S) \hat \otimes \ell_\infty(T,\ell_2)$.

  {\bf Note}. The spaces $\ell_\infty(S,\ell_2) \hat \otimes
\ell_\infty(T)$ and $\ell_\infty(S) \hat \otimes \ell_\infty(T,\ell_2)$ are
naturally continuously injected into $\ell_\infty(S \times  T,\ell_2)$, which
is used to define the above sum.

  {\bf Notation}. Let $\Omega = \tore^\nat$. Let $\mu$ be the normalized
Haar measure on $\Omega$, \ie $\mu = ({dt \over 2 \pi})^\nat$. We denote by $z
= (z_k)_{k \in \nat}$ a generic point of $\Omega$ (and we consider the $k-th$
coordinate $z_k$ as a function of $z$).

We will denote $\ell_\infty$ instead of $\ell_\infty(\nat)$ and
$\ell_\infty(\ell_2)$ instead of $\ell_\infty(\nat,\ell_2)$.

With this notation, we can state a more precise version of Theorem 2.1.

\proclaim Theorem 2.2. In the same situation as Theorem
2.1, let $\psi_1,...,\psi_n$ be a finite sequence in
$\ell_\infty(S) \hat \otimes \ell_\infty(T)$.\sn
(i) Assume
$$ \int \gamma_2(\sum_1^n z_k \psi_k) d\mu(z) < 1. \leqno (2.1)$$
Then there is a decomposition in $\ell_\infty(S) \hat \otimes \ell_\infty(T)$ of
the form
$$\psi_k = A_k + B_k$$
such that (with $\gamma$ as defined in (1.3))
$$\N {\sum_1^n e_k \otimes A_k}_{\ell_\infty(S,\ell_2) \hat \otimes_\gamma
\ell_\infty(T)} < 1$$
and
$$\N {\sum_1^n e_k \otimes B_k}_{\ell_\infty(S) \hat \otimes_\gamma
\ell_\infty(T,\ell_2)} < 1$$
(ii) Assume
$$\int \N {\sum_1^n z_k \psi_k}_{\ell_\infty(S) \hat \otimes \ell_\infty(T)}
d\mu(z) < 1$$
then there is a decomposition $\psi_k = A_k + B_k$ such that
$$\N {\sum e_k \otimes A_k}_{\ell_\infty(S,\ell_2)\hat \otimes \ell_\infty(T)}
< K_G~{\rm and}~\N {\sum e_k \otimes B_k}_{\ell_\infty(S) \hat \otimes
\ell_\infty(T,\ell_2)} < K_G,$$
where $K_G$ is the Grothendieck constant.

\proof The proof is based on the main result of [LPP] reformulated above as
Theorem 1.6.
By a standard Banach space technique, Theorem 2.2 can be reduced to the case when
$S$ and $T$ are finite sets. (Use the fact that $\ell_\infty$ is a ${\cal
L}_\infty$ space, more precisely it can be viewed as the
closure of the union of an increasing family of finite
dimensional sublattices each isometric to $\ell_\infty(S)$
for some finite set $S$).

We will denote by $\alpha_1$ (resp. $\alpha_2$) the norm on $\ell_1(S,\ell_2^n)
\otimes \ell_1(T)$ (resp. $\ell_1(S) \otimes \ell_1(T,\ell_2^n)$) which is dual
to the norm in $\ell_\infty (S,\ell_2^n) \hat \otimes_\gamma \ell_\infty(T)$
(resp. $\ell_\infty(S) \hat \otimes_\gamma \ell_\infty(T,\ell_2^n)$).

Let $(e_k)$ be the canonical basis of $\ell_2^n$.
Let $A_k \in \ell_1(S) \otimes \ell_1(T)$ and let $\Phi = \sum e_k \otimes A_k$.
 We will make the obvious
identifications permitting to view $\Phi$ as an element either of
$\ell_1(S,\ell_2^n) \otimes \ell_1(T)$ or of $\ell_1(S) \otimes
\ell_1(T,\ell_2^n)$.
Then, by duality Theorem 2.2 (i) is equivalent to the following inequality.

For
all $\psi_k$ in $\ell_\infty(S) \hat \otimes_{\gamma_2} \ell_\infty (T)$
$$|\sum < A_k,\psi_k>| \leq \int \gamma_2 (\sum z_k \psi_k) d\mu (z)
[\alpha_1(\Phi) + \alpha_2(\Phi)].\leqno (2.2)$$

To check this, by homogeneity we may assume (2.1) and also  
$$\alpha_1(\Phi) + \alpha_2(\Phi) = 1.\leqno (2.3)$$

Then, by Lemma 1.7, there are probabilities $P_1,P_2$ on $S$ and $Q_1,Q_2$ on
$T$ such that (with obvious identifications).
$$\forall~\ph \in \ell_\infty(S,\ell_2^n)~~~~~\forall~\eta \in
\ell_\infty(T)$$
$$|< \ph \otimes \eta,\Phi>| \leq \alpha_1(\Phi) \left(\int \N
{\ph(s)}_{\ell_2^n}^2 dP_1(s) \int |\eta(t)|^2 dQ_1(t)\right)^{1/2}$$
and
$$\forall~\beta \in \ell_\infty(S)~~~~\forall \omega \in
\ell_\infty(T,\ell_2^n)$$
$$|< \beta \otimes \omega,\Phi>| \leq \alpha_2(\Phi) \left(\int|\beta(s)|^2
dP_2(s) \int \N {\omega(t)}_{\ell_2^n}^2 dQ_2(t)\right)^{1/2}.$$

Now let $P = \alpha_1(\Phi) P_1 + \alpha_2(\Phi) P_2, Q = \alpha_1(\Phi) Q_1 + \alpha_2(\Phi) Q_2$. By
(2.3) these are probabilities.

Then
$$|< \ph \otimes \eta,\Phi>| \leq \left( \int \N {\ph(s)}_{\ell_2^n}^2
dP(s)\right)^{1/2} \left(\int|\eta(t)|^2 dQ(t)\right)^{1/2} \leqno (2.4)$$
and
$$|< \beta \otimes \omega,\Phi>| \leq \left(\int|\beta(s)|^2 dP(s)\right)^{1/2}
\left(\int\N {\omega(t)}_{\ell_2^n}^2 dQ(t)\right)^{1/2}.\leqno (2.5)$$

This means that $A_k$ defines a bounded linear operator $a_k :
L_2(P) \to L_2(Q)^*$ such that $< a_k(\beta),\eta> = < \beta
\otimes \eta,A_k>$. Moreover (2.5) and (2.4) imply
respectively $\N {\sum a_k^* a_k} \leq 1$ and $\N {\sum a_k
a_k^*} \leq 1$. Let $J_P : \ell_\infty(S) \to L_2(P)$ and
$J_Q : \ell_\infty(T) \to L_2(Q)$ be the canonical
inclusions. 
Let $\xi_k = (J_P \otimes J_Q) (\psi_k) \in L_2(P) \otimes
L_2(Q)$. By (1.9) we have
$$\int \N {\sum z_k \xi_k}_{L_2(P) \hat \otimes L_2(Q)} d\mu (z) < 1.$$
Note that $< \xi_k,a_k> = < \psi_k,A_k>$.
Hence applying (1.7) we obtain the desired inequality (2.2).
This concludes the proof of the first part. The second
part is an immediate consequence of the first one by
(1.6). \qed

 \n {\bf Remark} : It is also possible to deduce Theorem 2.2 directly from
the factorization Theorem of [P2] (see Corollary 1.7 or
Theorem 2.3 in [P2]), which applies in particular to
functions in $H^1$ with values in $\ell_\infty(S) \hat
\otimes_{\gamma_2} \ell_\infty(T)$. Using this, the
argument of [LPP] then gives the decomposition of Theorem
2.2 in a somewhat more explicit fashion as a formula in
terms of the factorization of the ``analytic" function $z
\to \sum z_k \psi_k$.

  \n{\bf Remark 2.3}. Let $A_k$ be as above such that 
$$\N {\sum_1^n e_k \otimes A_k}_{\ell_\infty(S,\ell_2) \hat \otimes_\gamma
\ell_\infty(T)} < 1.\leqno (2.9)$$
Then for any $n$-tuple $t_1,...,t_n$ in $T$ we have 
$$\sup_{s \in S} \left(\sum_{k=1}^n |A_k(s,t_k)|^2\right)^{1/2} < 1. \leqno
(2.10)$$
A similar remark holds for $\sum_1^n e_k \otimes B_k$.

Indeed, by the definition (1.3), (2.9) means that there is a Hilbert space $H$
and elements $\alpha$ in $\ell_\infty(S,\ell_2(H))$ and $\beta$ in
$\ell_\infty(T,H)$ each with norm $< 1$ such that
$$\forall~k=1,...,n~~~~~A_k(s,t) = < \alpha_k(s),\beta(t) > \leqno (2.11)$$
(where $\alpha(s) \in \ell_2(H)$ and $\alpha_k(s)$ denotes the $k-th$ coordinate
of $\alpha(s)$). Then (2.10) is an immediate consequence of (2.11).

We now derive Theorem 0.1 from Theorems 2.1 and 2.2.

  {\bf Proof of Theorem 0.1}. 

(i) $\Rightarrow$ (ii) is trivial. Assume (ii). By Remark 1.2 for almost all
choices of signs $\eps$ in $\{-1,1\}^G$ the function $(s,t) \to \eps(st)
\ph(st)$ defines a c.b. Schur multiplier $M_{\eps\ph}$ of
$B(\ell_2(G),\ell_2(G))$.

We can assume $\N {M_{\eps\ph}}_{cb} \leq F(\eps)$ for some measurable function
$F(\eps)$ finite almost everywhere on $\{-1,1\}^G$. A fortiori for each finite
subsets $S \subset G$ and $T \subset G$, the function $(s,t) \to \eps(st)
\ph(st)$ restricted to $S \times  T$ is a c.b. Schur multiplier of
$B(\ell_2(S),\ell_2(T))$ with norm $\leq F(\eps)$. By Proposition 1.1, this
means that we have for all finite subsets $S \subset G,T \subset G$
$$\N {\sum_{(s,t)\in S \times T} \eps(st) \ph(st) e_s \otimes
e_t}_{\ell_\infty(S) \hat \otimes_{\gamma_2} \ell_\infty(T)} \leq F(\eps),$$
where we have denoted by $(e_s)$ and $(e_t)$ the canonical bases of
$\ell_\infty(S)$ and $\ell_\infty(T)$.

Equivalently we have
$$\N {\sum_{x \in ST} \eps(x) \ph(x) \sum_{\buildrel{(s,t)\in S\times
T}\over{st=x}} e_s \otimes e_t}_{\ell_\infty(S) \hat
\otimes_{\gamma_2}\ell_\infty(T)} \leq F(\eps).\leqno (2.12)$$
By a classical integrability result of Kahane (\cf [Ka]) we
can assume that $F$ is integrable over $\{-1,1\}^G$, so
that there is a number  $C>0$ such that the average over
$\epsilon$ of the left side of (2.12) is less than $C$. By
a simple elementary reasoning (decompose into real and
imaginary parts, use the triangle inequality and the
unconditionality of the average over $\epsilon$), it
follows from (2.12) that if $\mu_G$ denotes the normalized
Haar measure on $\tore^G$ and if $z = (z_x)_{x \in G}$
denotes a generic point of $\tore^G$, we have $$\int \N
{\sum_{x\in ST}z_x \ph(x) \sum_{\buildrel{(s,t)\in S
\times T}\over {st=x}} e_s \otimes e_t}_{\ell_\infty(S) \hat
\otimes_{\gamma_2}\ell_\infty(T)} d\mu(z) < 2 C.$$ Let
$\psi_x = \ph(x) \sum_{\buildrel{(s,t)\in S\times
T}\over{st=x}} e_s \otimes e_t$ for all $x$ in $ST$ and let
$\psi_x=0$ otherwise.

Then by Theorem 2.2 and Remark 2.3 we have a decomposition 
$$\psi_x = A_x + B_x~{\rm in}~\ell_\infty(S) \otimes \ell_\infty(T)$$
such that
$$\sup_{s\in S} \left(\sum_{t \in T} |A_{st} (s,t)|^2 \right)^{1/2} < 2 C$$
and
$$\sup_{t\in T} \left(\sum_{s \in S} |B_{st} (s,t)|^2 \right)^{1/2} < 2 C$$
This yields functions $\ph_1$ and $\ph_2$ on $S \times  T$ such that
$\ph_1(s,t)=A_{st} (s,t),\ph_2(s,t)=B_{st} (s,t)$ and
$$\forall~(s,t) \in S \times  T~~~~~~\ph(st) = \psi_{st} (s,t) =
\ph_1(s,t) + \ph_2(s,t).$$
Hence
$$\sup_{s \in S} \left(\sum_{t \in T} |\ph_1(s,t)|^2 \right)^{1/2} < 2
C~~\sup_{t \in T} \left(\sum |\ph_2(s,t)|^2\right)^{1/2} < 2C.$$

Let us denote by $\ph_1^{S,T}$ and $\ph_2^{S,T}$ the functions obtained on $G
\times  G$ by extending $\ph_1$ and $\ph_2$ by zero outside $S \times  T$.

Now if
we let $S \times  T$ tend to $G \times  G$ along the set of all products of
finite sets directed by inclusion and if we let $\Phi_1,\Phi_2$ be pointwise
cluster points of the corresponding sets $(\ph_1^{ST})$ and $(\ph_2^{S,T})$, we
obtain finally two functions $\Phi_1$ and $\Phi_2$ on $G \times  G$ such that
$\ph(st) = \Phi_1(s,t) + \Phi_2(s,t)$ for all $(s,t)$ in $G \times  G$ and
satisfying 
$$\sup_{s \in G} \left(\sum_{t \in G} |\Phi_1(s,t)|^2\right)^{1/2} \leq 2
C~~~\sup_{t \in G} \left(\sum_{s \in G} |\Phi_2(s,t)|^2\right)^{1/2} \leq 2C$$
Let then $\Gamma_1 \cup \Gamma_2 = G \times  G$ be a partition defined by 
$$\Gamma_1 = \{(s,t) \in G \times  G| ~~~|\Phi_1(s,t)| \geq |\Phi_2(s,t)|\}$$
$$\Gamma_2 = \{(s,t) \in G \times  G| ~~~|\Phi_1(s,t)| < |\Phi_2(s,t)|\}$$
It is then clear that $\ph$ satisfies the property (iii) in Theorem 0.1. This
shows (ii) $\Rightarrow$ (iii). The equivalence (iii)
$\Leftrightarrow$ (iv) is   part
of Theorem 0.2 (due to Varopoulos). (Note that the
implication (iii) $\Rightarrow$ (i) also follows from the
implication (iii) $\Rightarrow$ (i) in Varopoulos's Theorem
0.2.) 

We now show (iii) $\Rightarrow$ (v). Assume (iii). Let
$a(x)$ be as in (v) and let $g$ and $h$ be in the unit
ball of $\ell_2(G,H)$. Assume
$$\max \left\{\N {\left(\sum
a(x)^* a(x)\right)^{1/2}},\N {\left(\sum a(x)
a(x)^*\right)^{1/2}}\right\}\leq 1.\leqno(2.14)$$It clearly
suffices to show that 
 $$|\sum_{s,t\in G\times G} \ph(st^{-1})<h(s),a(
st^{-1})g(t)>| \leq 2C.\leqno(2.15)$$
 Let $\Sigma_1$ (resp. $\Sigma_2$) be   the left side of
(2.14) with the summation restricted to $(s,t^{-1})\in\Gamma_1$
(resp. $(s,t^{-1})\in\Gamma_2$).
Observe that
$\sum_{(s,t^{-1})\in\Gamma_1}|\ph(st^{-1})|^2\|h(s)\|^2\leq
C^2$, hence by Cauchy-Schwarz and (2.14)   $$|\Sigma_1|\leq
C(\sum_{(s,t^{-1})\in \Gamma_1}  \|a(
st^{-1})g(t)\|^2)^{1/2} \leq C(\sum_{(s,t)\in
G\times G} <a(st^{-1})^* a(st^{-1})
g(t),g(t)>)^{1/2}$$
$$\leq C  (\sum_t \|g(t)\|^2 \|\sum_s
a(st^{-1})^* a(st^{-1})\|)^{1/2} \leq C.$$
 
A similar argument yields $|\Sigma_2|\leq C$ hence (2.15)
follows and the proof of (iii) $\Rightarrow$ (v) is
complete.

Finally we show (v) $\Rightarrow$
(i). We start by recalling that for any finitely supported
function $a : G \to B(H)$ we have the elementary inequality
$$\max \left\{\N {\left(\sum a(x)^* a(x)\right)^{1/2}},\N
{\left(\sum a(x) a(x)^*\right)^{1/2}}\right\} \leq \N {\sum
\lambda(x) \otimes a(x)}_{B(\ell_2(G,H))}.\leqno(2.16)$$
Now assume (v). We have then by (2.16) if $\sup \limits_x
|\eps(x)|\leq 1$ $$\N {\sum_{x \in G} \eps(x)\ph(x)
\lambda(x) \otimes a(x)} \leq C \N {\sum_{x \in G}
\lambda(x) \otimes a(x)},$$ hence the multiplier of
$C_\lambda^*(G)$ defined by $\eps\ph$ is completely bounded
with norm $\leq C$. This proves (v) $\Rightarrow$ (i). \qed

  {\bf Proof of Corollary 0.4}. Let $p$ and $\ph$ be as in Corollary 0.4.

Then there is a constant $C$ such that for all finite subsets $S,T$ of $G$ we
have
$$\int \N {\sum_{x \in p(S,T)} \ph(x) z_x (\sum_{p(s,t)=x} e_s \otimes
e_t)}_{\ell_\infty(G) \hat \otimes \ell_\infty(G)} < C.$$
Reasoning as above in the proof of (ii) $\Rightarrow$ (iii) in Theorem 0.1, we
find a decomposition of the form 
$$\ph(p(s,t)) = A_{p(s,t)} (s,t) + B_{p(s,t)} (s,t)$$
and using Remark 2.3 and the bounds
$$\sup_{s,x} |\{t |~~ p(s,t) = x \}|< \infty ~{\rm
and}~\sup_{t,x} |\{s|~~ p(s,t) = x \}|< \infty\leqno
(2.17)$$ we can obtain that for some constant $C'$
$$\sup_{s \in S} \sum_{t \in T} |A_{p(s,t)} (s,t)|^2 \leq C'$$
$$\sup_{t \in T} \sum_{s \in S} |A_{p(s,t)} (s,t)|^2 \leq C'.$$
We then conclude the proof as in the proof of Theorem 0.1 by a pointwise
compactness argument, showing that $(s,t) \to \ph(p(s,t))$ satisfies
(iii)' in Theorem 0.2. hence (recall Proposition 1.1 or
Theorem 1.3) for all bounded complex functions $(s,t) \to
\eps_{s,t}$ the function $(s,t) \to \eps_{s,t} \ph(p(s,t))$
is a Schur multiplier of $B(\ell_2(G),\ell_2(G))$.

\n{\bf Remark 2.4.}   Let
$S,T,X$ be arbitrary sets and let $p:S\times T\ra X$ be a
map satisfying (2.17). Consider a function  $\ph:X\ra \C$
and let $\psi:S\times T\ra \C$ be defined by
$$\psi(s,t)=\ph(p(s,t)).$$ Let $K(T)$ be the linear span
of the canonical basis of $\ell_2(T)$. For any $x$ in $X$,
let $\Lambda(x):K(T)\ra \C^S$ be the operator defined by
the matrix $\Lambda_x(s,t)$ defined by $\Lambda_x(s,t)=1$
if $p(s,t)=x$ and $\Lambda_x(s,t)=0$ otherwise.
Then we can generalize Varopoulos's theorem   as
follows. The following are equivalent\sn
(i) For all bounded functions $\eps : S\times T \to \comp$
the pointwise product $\eps \psi$ is a bounded Schur
multiplier of $B(\ell_2(S),\ell_2(T)$.\sn 
(ii) For almost
all choices of signs $\eps \in \{-1,1\}^{S\times T}$, the
product $\eps \psi$ is a bounded Schur
multiplier of $B(\ell_2(S),\ell_2(T)$.\sn
 (iii) There is
  a partition of $S \times  T$
 say $S \times  T = \Gamma_1 \cup \Gamma_2$ such that 
$$\sup_{s \in S} \sum_{t \in T} |\psi(s,t)|^2 1_{\{(s,t)
\in \Gamma_1\}} <\infty$$ $$\sup_{t \in T} \sum_{s \in S}
|\psi(s,t)|^2 1_{\{(s,t) \in \Gamma_2\}} <\infty.$$
\sn 
 (iv) There is
a constant $C$ such that for any Hilbert space $H$ and for
any finitely supported function $a : X \to B(H)$ we have
$$\N {\sum_{x \in X} \ph(x)\Lambda(x) \otimes
a(x)}_{B(\ell_2(S,H),\ell_2(T,H))} \leq C \max \left\{\N
{\left(\sum a(x)^* a(x)\right)^{1/2}},\N {\left(\sum a(x)
a(x)^*\right)^{1/2}}\right\}.$$
This statement is proved exactly as above. This applies
in particular when $p$ is the product map on a semigroup.
 When $S=T=X=\nat$ and $p(s,t)=s+t$ we recover results
already obtained in [B2].

  \vfill\eject
  {\bf \S~3. Lacunary sets.}

The study of ``thin sets" such as Sidon sets in discrete non Abelian groups has
been developped by several authors, namely Leinert [L1,L2] Bo$\dot{z}$ejko ([B1,
2, 3]), Figa-Talamanca and Picardello [FTP] and others. (See [LR] for the theory
of Sidon sets in Abelian groups).

In this section, we apply the preceding results to a class of ``lacunary"
sets which we call ``$L$-sets". There is a striking analogy between the
``$L$-sets" defined below and a class of subsets of $\nat$ which we will call
Paley sets. A subset $\Lambda \subset \nat$ will be called a Paley set if there
is a constant $C$ such that for all $f = \sum_{n=0}^\infty a_n e^{int}$ in
$H^1$ we have 
$$\left(\sum_{n \in \Lambda} |\hat f(n)|^2\right)^{1/2} \leq C \N
{f}_1.\leqno(3.1)$$ It is well known (\cf e.g. [R]) that
Paley sets are simply the finite unions of Hadamard-lacunary
sequences, \ie of sequences $\{n_k\}$ such that $\liminf_{k
\to \infty} (n_{k+1}/(n_k) > 1$.

Equivalently, $\Lambda$ is a Paley set iff there is a constant $C$ such that 
$$\forall~n > 0~~~~~~~|\Lambda \cap [n,2n[| \leq C.$$

In [LPP], it is proved that if $\Lambda \subset \nat$
is a Paley set there is a constant $C$ such that for any $f$ in $H^1(H \hat
\otimes H)$ there is a decomposition in $H \hat \otimes H$ of the form $\hat
f(n) = a(n) + b(n)~~~~\forall~n \in \Lambda$ such that
$$tr \left(\sum_{n \in \Lambda} a(n)^* a(n)\right)^{1/2} + tr \left(\sum_{n \in
\Lambda} b(n) b(n)^*\right)^{1/2} \leq C \N {f}_{H^1(H \hat \otimes
H)}\leqno(3.2)$$ Moreover, when $\hat f$ is supported by
$\Lambda$ this inequality becomes an equivalence.

The papers [LPP] and [HP] suggest that there is a strong analogy between Paley
sequences and free subsets of a discrete group $G$. To explain this we
introduce more notation. Let $H$ be a Hilbert space. We denote by $A(G,H \hat
\otimes H)$ the set of all functions $f : G \to H \hat \otimes H$ such that for
some $g,h$ in $\ell_2(G,H)$ we have
$$\forall~x \in G~~~~f(x) = \sum_{st=x} g(s) \otimes h(t).$$
Let
$$\N {f}_{A(G,H \hat \otimes H)} = \inf \{\N {g}_{\ell_2(G,H)} \N
{h}_{\ell_2(G,H)}\}$$
where the infimum runs over all possible representations. Then (see [HP]) if
$G$ is the free group on $n$ generators $g_1,...,g_n$, we have the following
analogue of (3.1). For any $f$ in $A(G,H \hat \otimes H)$
then there is a decomposition $f(g_k) = a_k + b_k$ in $H
\hat \otimes H$ such that  $$tr(\sum a^*_k a_k)^{1/2} + tr
(\sum b_k b^*_k)^{1/2} \leq 2 \N {f}_{A(G,H \hat \otimes
H)}.\leqno(3.3)$$
Moreover, when $  f$ is supported by
$\Lambda$ this inequality becomes an equivalence.

This   motivated   the following 

  {\bf Definition 3.1}. A subset $\Lambda$ of a discrete group $G$ will be
called an $L$-set if there is a constant $C$ such that for any $H$ and for any
$f$ in $A(G,H \hat \otimes H)$ we have
$$\inf_{f(x)=a(x)+b(x)} \left\{tr \left(\sum_{x \in \Lambda} a(x)^*
a(x)\right)^{1/2} + tr \left(\sum b(x) b(x)^*\right)^{1/2}\right\}$$
$$\leq C \N {f}_{A(G,H \hat \otimes H)},$$
where the infimum runs over all possible decompositions $f(x) = a(x) + b(x)$ in
$H \hat \otimes H$.

\proclaim Proposition 3.2. The following properties of a
subset $\Lambda \subset G$ are equivalent.\sn
(i) $\Lambda$ is an $L$-set.\sn
(ii) There is a constant $C$ such that for any Hilbert space $H$ and for any
finitely supported function $a : \Lambda \to B(H)$ we have
$$\N {\sum_{x \in \Lambda} \lambda(x) \otimes a(x)}_{B(\ell_2(G,H))} \leq C \max
\left\{\N {\left(\sum a(x)^* a(x)\right)^{1/2}},\N {\left(\sum a(x)
a(x)^*\right)^{1/2}}\right\}.$$\sn
(iii) For any bounded sequence $\eps$ in
$\ell_\infty(G)$ supported by $\Lambda$, the associated multiplier defined on
$C^*_\lambda(G)$ by  $$\sum f(x) \lambda(x) \to \sum f(x) \eps(x) \lambda(x)$$
is completely bounded on $C^*_\lambda(G)$.\sn
 
\proof (i) and (ii) are clearly equivalent. They are but a dual reformulation of
each other.
(ii) and (iii) are   equivalent by Theorem 0.1 applied to
the indicator function of $\Lambda$. \qed

  \n{\bf Remark}. By Theorem 0.1, the preceding properties are also equivalent
to the property (iii)'   obtained by requiring that
the property (iii) holds only for almost all choice of signs
$(\eps(x))_{x \in \Lambda}$ in $\{-1,1\}^\Lambda$.  

  \n{\bf Remark}. The preceding result shows that $\Lambda$ is
an $L$-set iff $\Lambda$ is a strong $2$-Leinert set in the
sense of Bo$\dot{z}$ejko [B1]. Leinert [L1,L2]    
 first
constructed infinite sets of this kind in free
noncommutative groups. Leinert's  results were clarified
in [AO]. Moreover, in [H1] several related important
inequalities were obtained for the operator norm  of the
convolution on the free group
  by a function supported by the words of a given fixed
lenghth in the generators. In particular, it was known to
Haagerup (see [HP]) that any free subset of a discrete
group is an $L$-set. For instance the generators
(or the words of
length one) on the
free group with countably many generators  form an $L$-set.
On the other hand it is rather easy to see
that the set of words of a fixed length $k>1$ is not an $L$-set
(for instance it clearly does not satisfy (ii) in Theorem 3.3).  

The main result of this section is the following 

\proclaim Theorem 3.3. Let $G$ be an arbitrary discrete
group. Let $\Lambda \subset G$ be a subset. Let $R_\Lambda
\subset G \times  G$ be defined by  $$R_\Lambda =\{(s,t)
\in G \times  G |\  st \in \Lambda\}.$$ The following
properties of $\Lambda$ are equivalent \sn (i) $\Lambda$ is an
$L$-set.\sn (ii) There is a constant $C$ such that for any
finite subsets $E,F \subset G$ with $|E| = |F| = N$ we have 
$$|R_\Lambda \cap (E \times  F)|\leq C N.$$\sn
(iii) There is a constant $C$ and there is a partition
$R_\Lambda = \Gamma_1 \cup \Gamma_2$ such that $$\sup_{s
\in G} \sum_{t \in G} 1_{(s,t) \in \Gamma_1  
 }\leq C\quad {\rm  and}\quad  \sup_{t \in G}\sum_{s \in G}
1_{(s,t) \in \Gamma_2 }\leq C.$$

\p
This is an immediate consequence of   Theorem 0.1 applied to the indicator
function of $\Lambda$.

 \n {\bf Remark}. As already mentioned, the equivalence between (ii) and
(iii) is due to Varopoulos [V1]. Our result shows that $\Lambda \subset G$ is
an $L$-set iff $R_\Lambda$ determines a $V$-set for $\ell_\infty(G) \hat \otimes
\ell_\infty(G)$ in the sense of Varopoulos (see [LP] and
[V2]). Equivalently, let $G' = \tore^G$ and let
$(\gamma_s)_{s \in G}$ be the coordinates on $G$. Then
$\Lambda$ is an $L$-set in $G$ iff the set $\{\gamma_s
\times  \gamma_t|st \in \Lambda\}$ is a Sidon set in the
dual of $G' \times  G'$ \ie in $\ent^{(G)} \times 
\ent^{(G)}$.

\n{\bf Remark}. Taking Remark 2.4 into account, we can
extend the notion of $L$-set to the case when $\Lambda$ is
a subset of a semi-group $G$ embeddable into a group (for
instance $G=\nat$). In that case we will say that  $\Lambda$ 
is an $L$-set if it satisfies the equivalent properties
of Theorem 3.3 (or the analogue of (ii) in Proposition
3.2). This provides a common framework for Paley sets and
$L$-sets.   Note however that all the   $L$-subsets of
$\ent$ (or of any amenable group) are finite, while the
$L$-subsets of $\nat$ are exactly the Paley sets. Thus this
notion   depends of the choice of the
semi-group containing $\Lambda$. If we remain in the
category of groups this difficulty does not arise, if $H$
is a subgroup of a group $G$ and if $\Lambda\subset H$,
then $\Lambda$ is an $L$-set in $H$ iff it is an $L$-set in
$G$ (this is easy to check e.g. by Proposition 3.2).

Note that $L$-sets are clearly stable under finite unions. Moreover the
translate of an $L$-set is again an $L$-set. The only known examples of $L$-sets
seem to be finite unions of translates of free sets. The sets which are
translates of a free set (more precisely translates of a free set augmented by
the unit element)  are characterized in [A0] as those which
have the Leinert property. In analogy with Paley sequences
we formulate the following.

    {\bf Conjecture 3.5}. Every $L$-set $\Lambda$ can be written as a finite
union $x_1 F_1 \cup...\cup x_n F_n$ where $x_1,...,x_n \in G$ and $F_1,...,F_n$
are free subsets of $G$.(Here,   the subset reduced to the unit
element is considered free, so that a singleton is a translate of a free set.)

It is possible to check that if $\Lambda$ satisfies (iii)
in Theorem 3.3 with $C=1$ then it satisfies the Leinert
property in the sense of [AO], hence it is a translate of
a free set augmented by the unit element, a fortiori it is
the union of two translates of free sets. Therefore to
verify the above conjecture it suffices to prove that any
set $\Lambda$ satisfying (iii)
in Theorem 3.3 with some constant $C$ can be written as a
finite union of sets satisfying the same property with
$C=1$.

\vfill\eject
  {\bf \S~4. A more general framework.}

Actually, Theorem 2.2 can be extended to a rather general situation already
considered in [P2]. We describe this briefly since it is easy to adapt the
preceding ideas to this setting.

Let $X$ be a Banach space.
We will identify an element $u$ in $X   \otimes \ell_2$ with an operator $u
: X^* \to \ell_2$ (of finite rank and weak-* continuous). Hence, for any $\xi$
in $X^*,u(\xi) \in \ell_2$.

A norm $\delta$ on $X \otimes \ell_2$ will 
be called 2-convex if there is a constant $c>0$
such that for any $u$ in $X   \otimes \ell_2$
$$c\|u\|=c\sup_{\xi\in X^*,\|\xi\| \leq 1}\|u(\xi)\| \leq
\delta(u)\leqno(4.1)$$ and such that   for all $u,u_1,u_2$
in $X \otimes \ell_2$ satisfying  $$\forall~\xi \in
X^*~~~~\N {u(\xi)} \leq \left(\N {u_1(\xi)}^2 + \N
{u_2(\xi)}^2\right)^{1/2},$$ we have
$$\delta(u) \leq (\delta(u_1)^2 + \delta \left(u_2)^2\right)^{1/2}.\leqno(4.2)$$
Note that if $\N {u(\xi)} = \N {u_1(\xi)}$ for all $\xi$ in $X^*$, we must have
$\delta(u) = \delta(u_1)$, moreover for all $T:\ell_2\ra
\ell_2$ we have $\delta(Tu)\leq \|T\|\delta (u).$

Now let $X,Y$ be two Banach spaces and let $\delta_1$ (resp. $\delta_2$) be a
2-convex norm on $X \otimes \ell_2$ (resp. $Y \otimes \ell_2$). We can
introduce a norm $\Gamma$ on $X \otimes Y$ by setting $\forall~u \in X \otimes
Y,u = \sum_{i=1}^n x_i \otimes y_i$
$$\Gamma(u) = \inf \left\{\delta_1(\sum_{i=1}^n x_i \otimes e_i) \delta_2
(\sum_{i=1}^n y_i \otimes e_i)\right\}$$
where the infimum runs over all possible decompositions of $u$ and where $e_i$
denotes the canonical basis of $\ell_2$. It is easy to see that this is a norm.
We denote by $X \hat \otimes_\Gamma Y$ the completion of $X \otimes Y$ for this
norm.

We also denote by $X \hat \otimes_{\delta_1} \ell_2$ and $Y \hat
\otimes_{\delta_2} \ell_2$ the completions of $X \otimes \ell_2$ and $Y \otimes
\ell_2$ for the norms $\delta_1$ and $\delta_2$.

Assume that $X$ and $Y$ are continuously injected in a Banach space $Z$.
Then, by (4.1)  $X \hat \otimes_{\delta_1} \ell_2$ and $Y
\hat \otimes_{\delta_2} \ell_2$ are both continuously
injected into $Z \buildrel{\vee} \over{\otimes}\ell_2$ (the
injective tensor product), so that we can give a meaning to
the sum $X \hat \otimes_{\delta_1} \ell_2 + Y \hat
\otimes_{\delta_2} \ell_2$.

For simplicity, we denote
$$X [\ell_2] = X \hat \otimes_{\delta_1} \ell_2~{\rm and} ~Y[\ell_2] = Y \hat
\otimes_{\delta_2}\ell_2.$$
We can now equip $X[\ell_2] \otimes \ell_2$ with a 2-convex norm $\Delta_1$ as
follows. For any $v = \sum_{i=1}^n v_i \otimes e_i$ with $v_i \in X[\ell_2]$
there is clearly an operator $w \in X \hat \otimes_{\delta_1} \ell_2$ (not
necessarily unique) such that $\N {w (\xi)} = \left(\sum_{i=1}^n \N
{v_i(\xi)}^2\right)^{1/2}~\forall~\xi \in X^*$. We then define
$$\Delta_1(v) = \delta_1(w).$$ By (4.2) this does not
depend on the particular choice of $w$. Clearly, this
defines (by the density of $\cup \ell_2^n$ in $\ell_2$) a
2-convex norm $\Delta_1$ on $X[\ell_2] \otimes \ell_2$.
Now using the pair $(\Delta_1,\delta_2)$ (instead of the
pair $(\delta_1,\delta_2))$ we can define the space
$$X[\ell_2] \hat \otimes_\Gamma Y$$ exactly as above for $X
\hat \otimes_\Gamma Y$.

Similarly, we can define $\Delta_2$ on $Y[\ell_2]\otimes \ell_2$ and using the
pair $(\delta_1,\Delta_2)$ we construct the space
$$X \hat \otimes_\Gamma Y [\ell_2]$$
exactly as above for $X \hat \otimes_\Gamma Y$.

Assume that $X,Y$ are both continuously injected in a Banach space $Z$.
Then, by (4.1),  it is easy to check that $X \hat
\otimes_{\delta_1} \ell_2$ and $Y \hat
\otimes_{\delta_2}\ell_2$ are both continuously injected
into the injective tensor products $X \buildrel{\vee}\over
{\otimes} \ell_2$ and $Y \buildrel {\vee}\over{\otimes}
\ell_2$, and consequently also $X[\ell_2] \hat
\otimes_\Gamma Y$ and $X \hat \otimes_\Gamma Y [\ell_2]$
are both continuously injected into $X \buildrel
{\vee}\over{\otimes} Y \buildrel {\vee}\over {\otimes}
\ell_2$, so that using this inclusion we may consider the
sum $$X [\ell_2] \hat \otimes_\Gamma Y + X \hat
\otimes_\Gamma Y [\ell_2],$$ with its natural norm (see
section 1).

Then Theorems 2.1 and 2.2 can be generalized as follows.

\proclaim Theorem 4.1. With the preceding notation, consider elements $\psi_n$
in $X \hat \otimes_\Gamma Y$ such that the series $\sum_{n=1}^\infty \eps_n
\psi_n$ converges for almost all choice of signs $\eps_n = \pm 1$. Then
necessarily the series $\sum_{n=1}^\infty e_n \otimes \psi_n$ converges in the
space $X [\ell_2] \hat \otimes_\Gamma Y + X \hat \otimes _\Gamma Y [\ell_2]$.

\proclaim Theorem 4.2. Let $\psi_1,...,\psi_N$ be a finite sequence in $X \hat
\otimes_\Gamma Y$ such that 
$$\int \Gamma (\sum_1^N z_k \psi_k) d\mu(z) < 1.$$
Then there is a decomposition in $X \hat \otimes_\Gamma Y$ of the form 
$$\psi_k = A_k + B_k$$
such that 
$$\N {\sum_1^N e_k \otimes A_k}_{X[\ell_2]\hat \otimes_\Gamma Y} < 1$$
and
$$\N {\sum_1^N e_k \otimes B_k}_{X \hat \otimes_\Gamma Y [\ell_2]} < 1.$$

(  {\bf Note}. Here of course $\sum e_k \otimes A_k$ is identified with
an element of $X[\ell_2]\hat \otimes_\Gamma Y$ in the obvious natural way and
similarly for $\sum e_k \otimes B_k$.)

The proof is the same as for Theorems 2.1 and 2.2. We leave the details to the
reader.

 \n {\bf Remark 4.3}. In particular, with the notation and terminology of
[LPP]we find (without any UMD assumption) that if $X,Y$ are two 2-convex
Banach lattices, then there is a natural inclusion Rad$(X \hat \otimes Y)
\subset X(\ell_2) \hat \otimes Y + X \hat \otimes Y(\ell_2)$, so that if $X,Y$
are both of finite cotype we have 
$${\rm Rad}(X \hat \otimes Y) \approx ~{\rm Rad}(X) \hat \otimes Y + X \hat
\otimes ~{\rm Rad}(Y).$$
This is a slight refinement of some of the results of [LPP].

\baselineskip=18pt
\def\qed{{\vrule height7pt width7pt
depth0pt}\par\bigskip}
\def\n{\noindent}
\def\Z{\ent}
\def\T{\tore}
\def\p{{\bf Proof:\ }}
\def\C{\comp}
\def\ra{\rightarrow}
\def\p.{{\bf Proof:\ }}
 
\n {\bf Remark 4.4.} As in Remark 1.9, let us denote
 by $H_1(\T; X)$ (resp. $H_1({\T}^{\nat};X)$) the
subspace of   $L_1(\T,m;X)$ (resp. $L_1({\T}^{\nat},{m}^{\nat};X)$) formed by
 the functions with Fourier transform supported by
the non-negative elements. We define similarly (for short)
the space $H_{\infty}(\T; X)$ (resp.
$H_{\infty}({\T}^{\nat};X)$). Using
Remark 1.9 we obtain the same conclusion as Theorem 4.2 whenever there is a function $f$ in the
interior of the  unit ball of $H_1(\T;X\otimes_\Gamma Y)$   (resp. $H_1({\T}^{\nat};X\otimes_\Gamma
Y)$) such that $\hat{f}(3^k)= \psi_k$ (resp. $\hat{f}(z_k)= \psi_k$) for all $k=1,...,n$.
 Of course this remark applies in particular to Theorem 2.2. In the case of
Theorem 2.2 this remark seems useful because it turns out the converse is true.
More precisely using a classical inequality (cf. [R] p.222) it can be proved
that, in the situation of Theorem 2.2, for every element in the unit ball of 
${\ell_\infty(S,\ell_2)\hat \otimes \ell_\infty(T)}$ with
a finitely supported sequence of coefficients $(A_k)$ 
in ${\ell_\infty(S )\hat \otimes \ell_\infty(T)}$ there is a function
$f$ in $H_{\infty}(\T;{\ell_\infty(S )\hat
 \otimes \ell_\infty(T)})\subset H_{1}(\T;{\ell_\infty(S
)\hat
 \otimes \ell_\infty(T)})$ with norm
less than an absolute constant $C$ such that 
$$\hat{f}(3^k)= A_k.$$ We can treat similarly any element
in the unit ball of 
${\ell_\infty(S)\hat \otimes \ell_\infty(T,\ell_2)}$,
hence the same conclusion holds for any element (with
a finitely supported sequence of coefficients) in the
unit ball of the space $${\ell_\infty(S,\ell_2)\hat \otimes
\ell_\infty(T)}+{\ell_\infty(S)\hat \otimes
\ell_\infty(T,\ell_2)}.$$ In other words, the point of the
present remark is that it
 yields a characterization of the sequences $(\psi_k)$ for which the conclusion
of Theorem 2.2 holds.

 \vfill \eject

 \baselineskip=18pt
\def\qed{{\vrule height7pt width7pt
\depth0pt}\par\bigskip}

\def\T{{\bf T}}
\def\n{\noindent}
\def\p{{\bf Proof:\ }}
\def\C{\comp}
\def\ra{\rightarrow}
\def\p.{{\bf Proof:\ }}
\S~{\bf 5. More applications to completely bounded maps}.

As emphasized in [BP] (see also [P4]) the $cb$ norm on $B(M_n,B(H))$ can be
viewed as an example of the $\Gamma$ norms discussed in section 4. In particular
we can obtain an analogue of Theorem 2.2 for c.b. maps.

\proclaim Theorem 5.1. Consider Hilbert spaces $H_1$ and $H$ and completely
bounded maps
$$u_1,...,u_N : B(H_1) \to B(H).$$
Assume that $\int \N {\sum \limits_{k=1}^N z_k u_k}_{cb} d\mu(z) < 1$. Then there
is for some Hilbert space $K$ a representation 
$$\pi : B(H_1) \to B(K)$$
and operators $V_k : H \to K,W_k : H \to K,V : H \to K~~W : H \to K$ such that
$$\N {V} \leq 1,\N {W} \leq 1,\N {\sum_1^N {V_k}^* V_k} \leq 1~~~\N {\sum_1^N
{W_k}^* W_k} \leq 1$$
and such that for all $k = 1,...,N$
$$\forall ~x \in B(H_1) ~~u_k(x) = V_k^* \pi(x) W + V^*\pi(x) W_k.$$

\n{\bf Remark 5.2}.
First observe that it is enough to find representations $\pi' : B(H_1) \to
B(K')~~~\pi'' : B(H_1) \to B(K'')$ such that $u_k(.) = V_k^* \pi'(.)W + V^* 
\pi''(.) W_k$ since we can replace each of $\pi'$ and $\pi''$ by $\pi' \oplus
\pi''$.

\n{\bf Remark 5.3}.
Assume that we have a net $(u^\alpha_k)_{k \leq N}$ of $N$-tuples of maps from
$B(H_1)$ into $B(H)$ such that for all $x$ in $B(H_1)~~~u^\alpha_m(x) \to
u_k(x)$ when $\alpha \to \infty$ and satisfying the conclusions of Theorem 5.1,
\ie such that there is a Hilbert space $K_\alpha$, a representation $\pi_\alpha
: B(H_1) \to B(K_\alpha)$ and operators
$V^\alpha_k,V^\alpha,W^\alpha_k,W^\alpha$such that
$V^\alpha,W^\alpha,\sum \limits_k V^{\alpha*}_k V^\alpha_k$ and $\sum \limits_k
W^{\alpha*}_k W^\alpha_k$ are all of norm $\leq 1$ and we have for all $x$ in
$B(H_1)$
$$u^\alpha_k(x) = V^{\alpha *}_k \pi_\alpha(x) W_k + V^{\alpha *} \pi_\alpha(x)
W^\alpha_k.$$
Then $\{u_k\}$ satisfies the conclusions of Theorem 5.1. Indeed, we can take
for $K$ an ultraproduct of the Hilbert spaces $K_\alpha$, and similarly for
$\pi$ and  for the operators $V,W,V_k,W_k$.

\n{\bf Remark 5.4}.
Assume $H_1$ and $H$ both finite dimensional. Then any operator $u : B(H_1) \to
B(H)$ can be identified with a linear operator $\tilde u : H_1 \otimes H \to
(H_1 \otimes H)^*$ defined by 
$$\forall~x_1,y_1 \in H_1~~~~~\forall~x,y \in H$$
$$< \tilde u(x_1 \otimes x),y_1 \otimes y > =  < \overline {y},u(x_1 \otimes
y_1)x >.$$
({\bf Note}: We denote by $y \to \overline {y}$ an anti isometry of $H$
onto itself ; note that on the left side we have a bilinear pairing while the
scalar product appearing on the right side is antilinear in the first variable).

Consider a factorization of $\tilde u$ of the form 
$$\tilde u = \sum_{i=1}^{i=n} x_i \otimes y_i$$
with $x_i,y_i \in (H_1 \otimes H)^*$

We define
$$\delta_1(\sum_{i=1}^n x_i \otimes e_i) = \sup \{|\sum_{i=1}^n < x_i,h_i
\otimes h>| h_i \in H_1, \sum \N {h_i}^2 \leq 1, h \in H,\N {h} \leq
1\}.\leqno (5.1)
$$ We may identify an element $x_i$ in $(H_1 \otimes H)^*$ with a
linear operator $V_i : H \to H_1$ by setting
$$\forall~k \in H,\forall~h \in H~~~<x_i,k \otimes h > = < \overline {k},V_i h
>\leqno (5.2)$$
Then (5.1) becomes
$$\eqalign{
\delta_1(\sum_{i=1}^n x_i \otimes e_i)& = \sup \{(\sum \N {V_i h}^2)^{1/2}|h
\in H~~~\N {h} \leq 1\}\cr
& = \N {(\sum V_i^* V_i)^{1/2}}}{.}\leqno (5.3)$$
Let $X = (H_1 \otimes H)^*$ and let ${\cal V}$ be the linear span of the basis
vectors of $\ell_2$. Clearly the formula (5.3) defines a 2-convex norm on $X
\otimes \ell_2$ (by density, say, of $X \otimes {\cal V}$ in $X \otimes \ell_2)$.

We set $Y = X$ and we define $\delta_2 = \delta_1$ on $Y \otimes \ell_2$. Then we
can define the norm $\Gamma$ associated to $\delta_1$ and $\delta_2$ as in
section 4, and also the norms $\Delta_1$ and $\Delta_2$ and the spaces $X
[\ell_2] \otimes_\Gamma Y$ and $X[\ell_2]\otimes_\Gamma Y$.

By well known results on the factorization of $c.b.$ maps (\cf [Pa]) we have then
$$\N {u}_{cb} = \Gamma(\tilde u).$$
Moreover, if $u_1,...,u_N$ are given $cb$ maps from $B(H_1)$ into $B(H)$, and
if $\tilde u_1,...,\tilde u_N$ denote the corresponding elements of $X \otimes
Y$ (with $X = Y = (H \otimes H_1)^*).$ We claim that 
$$\N {\sum_{k=1}^N \tilde u_k \otimes e_k}_{X[\ell_2]\otimes_\Gamma Y} <
1\leqno (5.4)$$
iff there are operators $\{W_i^k |k \leq N,i \leq n\}$ and $\{V_i |i \leq n\}$
such that $\forall~x \in B(H_1)$ 
$$u_k(x) = \sum_{i=1}^n V_i^* x W_i^k~ {\rm and}~\N {\sum_i V_i^* V_i} < 1 \N
{\sum_{i,k} {W_i^{k}}^* W_i^k} < 1.$$
Indeed (5.4) holds iff we can write
$$\sum_{k=1}^N \tilde u_k \otimes e_k = \sum_{i=1}^n \xi_i \otimes \eta_i$$
with $\xi_i \in X \otimes \ell_2^N$ and $\eta_i \in Y$ such that
$$\Delta_1(\sum \xi_i \otimes e_i) < 1$$
and
$$\delta_2(\sum \eta_i \otimes e_i) < 1.$$
Let $\xi_i = \sum \limits_{k \leq N} \xi_i^k \otimes e_k$ with $\xi_i^k \in X$
and let $V_i^k$ and $W_i$ be the operators associated to $\xi_i^k$ and $\eta_i$
by the correspondence (5.2). We then obtain the above claim.

This remark shows that in the case when both $H_1$ and $H$ are finite
dimensional, Theorem 5.1 can be viewed as a particular case of Theorem 4.1.

\n{\bf Proof of Theorem 5.1}.
We assume $H$ finite dimensional until the last step of the proof. The case
when $H_1$ is also finite dimensional has been checked in the preceding remark.
Assume $H_1$ infinite dimensional assume that each $u_1,...u_N$ is
weak*-continuous, \ie continuous from $\sigma(B(H_1),B(H_1)_*)$ into
$B(H)$. We may use the fact that there is an increasing set $B(H_\alpha)
\subset B(H_1)$ with $\dim H_\alpha < \infty$ such that $\bigcup \limits_\alpha
B(H_\alpha)$ is weak*-dense in $B(H_1)$. Applying the first part of the proof
to the restrictions ${{u_k}_{|B(H_\alpha)}}$ for each $\alpha$ and passing to
the limit in a standard way (as in remark 5.3) we obtain Theorem 5.1 in that
case also.

Next when $H_1$ is arbitrary and $H$ finite dimensional we can involve the local
reflexivity principle (\cf e.g. [D]) to claim that there is a net
$(u_k^\alpha)_{k \leq N}$ of maps which are weak* continuous from $B(H_1)$ into
$B(H)$, which tend pointwise to $(u_k)_{k \leq N}$ when $\alpha \to \infty$ and
which satisfy
$$\int \N {\sum z_k u_k^\alpha}_{cb} d\mu(z) < 1.$$
By remark 5.3 we obtain Theorem 5.1 in that case also.

Finally we remove the assumption that $H$ is finite dimensional.

Let $H_\alpha$ be an increasing net of finite dimensional subspaces of $H$ with
$\overline {\bigcup H_\alpha} = H$.

For each $k$ and $\alpha$ let
$$u_k^\alpha(x) = P_{H_\alpha}u_{k}(x)_{|H_\alpha} \in B(H_\alpha).$$
By the first part of the proof the conclusion of Theorem 5.1 holds for
$(u_k^\alpha)_{k \leq N}$ for each $\alpha$.

Using an ultraproduct argument as in Remark 5.3 we conclude one more that
$u_1,...,u_N$ satisfy the conclusions of Theorem 5.1.

We now give some consequences of Theorem 5.1. We denote again by $\ell_2^n$ the
$n$ dimensional Hilbert space equipped with its canonical basis $(e_i)_{i \leq
n}$. We will identify $M_n$ with $B(\ell_2^n)$ and $e_{ij}$ with $e_i \otimes
e_j$.

We will need the following notation.

\n Let $H_1,H_2$ be two Hilbert spaces and let
$S_1 \subset B(H_1),S_2 \subset B(H_2)$ be closed subspaces (``operator
spaces"). We will denote by $S_1 \otimes_{\min} S_2$ the completion of $S_1
\otimes S_2$ equipped with the norm induced by $B(H_1 \otimes_2 H_2)$, where
$H_1 \otimes_2 H_2$ is the Hilbertian tensor product. The space $S_1
\otimes_{\min} S_2$ is called the minimal (or the spatial) tensor
product of $S_1$ and $S_2$. In particular we have obviously
$M_n(B(H))=M_n\otimes_{\min}B(H)$.

 \n We will denote by $BR_n$ (resp.$BC_n$) the subspace of $M_n \otimes_{\min} M_n$
formed by all elements of the form
$$~~~~~~~~~~~y_1 = \sum_i e_{ii}\otimes \sum_j x_{ij} e_{ij})$$
$${\rm (resp}.~~~~~~~ y_2 = \sum_j e_{jj} \otimes \sum_i x_{ij} e_{ij}.$$
Note that $\N {y_1} = \sup \limits_j(\sum \limits_i |x_{ij}|^2)^{1/2}$ and $\N
{y_2} = \sup \limits_i (\sum \limits_j |x_{ij}|^2)^{1/2}$ so that $BR_n$ (resp.
$BC_n$) is naturally isometric with the space of matrices with ``bounded rows"
(resp. ``bounded columns"), which explains our notation.

\n Let $J_1 : M_n \to BR_n$ (resp. $J_2 : M_n \to BC_n)$ be the map defined by
$$~~~~~~~~~J_1(x) = \sum_i e_{ii} \otimes \sum_j x_{ij}
e_{ij}$$

$${\rm(resp.}~~~~~J_2(x) = \sum_j e_{jj} \otimes \sum_i x_{ij} e_{ij}).)$$
Obviously we have $\N {J_1} \leq 1$ and $\N {J_2}\leq 1$. Moreover a simple
verification shows that 
$$\N {J_1}_{cb} \leq 1~{\rm and}~\N {J_2}_{cb} \leq 1.\leqno (5.5)$$
Let us denote here $G = {\bf T}^{n^2}$ and let $\mu$ be the normalized Haar
measure on $G$. By a simple computation one can chech that for any Hilbert space
$H$ and for any $x_{ij}$ in $B(H)$ we have 
$$\eqalign{\N {\sum_i e_{ii} \otimes \sum_j e_{ij} \otimes x_{ij}}_{BR_n\otimes_{\min} B(H)}
&=
\sup_i \N {\sum_j e_{ij} \otimes x_{ij}}_{M_n( B(H))}\cr
& = \sup_i \N {(\sum_j
x_{ij}x_{ij}^*)^{1/2}}_{B(H)}}\leqno(5.6)'$$
and similarly
$$\N {\sum_j e_{jj}\otimes \sum_i e_{ij}\otimes x_{ij}}_{BC_n\otimes_{\min} B(H)}
= \sup_j \N {(\sum_i x_{ij}^* x_{ij})^{1/2}}_{B(H)}.\leqno (5.6)''$$
Observe that the preceding expressions do not change if we replace $(x_{ij})$
by $(z_{ij} x_{ij})$ with $|z_{ij}| = 1$, i.e. with $z=(z_{ij})\in G$. Hence if
we denote by $T_z : M_n \to M_n$ the Schur multiplier defined by 
$$T_z((a_{ij})) = (a_{ij}z_{ij}),$$
we find using (5.5) and this observation that for all $z = (z_{ij})$ in $G$ we
have
$$\N {J_1 T_z}_{cb} \leq 1~~~~~~~\N {J_2 T_z} \leq 1.\leqno(5.7)$$
We can now state

\proclaim Corollary 5.5. Let $H$ be a Hilbert space. Consider an operator $u :
M_n \to B(H)$. Let $u_z = u T_z : M_n \to B(H)$. 
\item(i) Assume $\int \N {u_z}_{cb}
d\mu(z) < 1$. Then there are operators $a_1 : BR_n \to B(H)$ and $a_2 : BC_n \to
B(H)$ such that   $$\N {a_1}_{cb} \leq 1~~~~~\N
{a_2}_{cb} \leq 1 \ {\rm and}\ u = a_1 J_1 + a_2 J_2.\leqno (5.8)$$ 
\item(ii) Conversely if (5.8) holds we have
$$\int \N {u_z}_{cb} d \mu(z) \leq \sup_{z \in G} \N {u_z}_{cb} \leq 2. \leqno
(5.9)$$

\n{\bf Proof}. Note that (5.9) is an obvious consequence of (5.7) so it suffices to
prove the first part. Let $u_{ij}: M_n \to B(H)$ be defined by $u_{ij}(x) =
x_{ij}u(e_{ij})$ so that $u_z = \sum \limits_{ij}z_{ij}u_{ij}$. By Theorem
5.1 we can find a Hilbert space $K$ operators $V^{ij},W^{ij},V,W$ from $H$ into
$K$ and a representation $\pi : M_n \to B(K)$ such that for all $x$ in $M_n$
$$u_{ij}(x) = V^{ij*} \pi(x) W + V^* \pi(x) W^{ij} \leqno (5.10)$$
and
$$\N {V} \leq 1,\N {W}\leq 1,\N {\sum_{ij} V^{ij*} V^{ij}} \leq 1,\N {\sum_{ij} W^{ij*}
W^{ij}}\leq 1.$$
Let $T_{ij} = u(e_{ij})$. We deduce from (5.10)
$$T_{ij}= V^{ij*} \pi(e_{ij}) W + V^* \pi(e_{ij}) W^{ij}.\leqno (5.11)$$
Let $a_1 : BR_n \to B(H)$ and $a_2 : BC_n \to B(H)$ be defined by 
$$~~~~~~~~a_1(e_{ii}\otimes e_{ij}) = V^* \pi(e_{ij}) W^{ij}~~~~~{\rm and}$$
$$a_2(e_{jj} \otimes e_{ij}) = V^{ij*}\pi(e_{ij}) W.$$
Clearly by (5.11) we have $u = a_1 J_1 + a_2 J_2$. We claim that $\N {a_1}_{cb}
\leq 1$ and $\N {a_2}_{cb} \leq 1$. To check this we will use the well known
inequality 
$$\N {\sum a_k^* b_k} \leq \N {(\sum a_k^* a_k)^{1/2}} \N {(\sum b_k^*
b_k)^{1/2})} \leqno(5.12)$$
valid for $a_k,b_k \in B(H)$.
For any $\xi_{ij}$ in $B(K_1)~~~(K_1$ an arbitrary Hilbert space) we have
$$\eqalign{&\N {\sum a_1(e_{ii} \otimes e_{ij}) \otimes \xi_{ij}}= \N {\sum V^*
\pi(e_{ij}) W^{ij}\otimes \xi_{ij}}\cr
\leq &\N {(\sum W^{ij*} W^{ij})^{1/2}}\N {(\sum (V^* \pi (e_{ij}) \otimes
\xi_{ij}) (V^* \pi (e_{ij}) \otimes \xi_{ij})^*)^{1/2}}\cr
\leq & \N {V} \N {(\sum_{ij}\pi (e_{ij})\pi (e_{ij})^* \otimes
\xi_{ij}\xi_{ij}^*)^{1/2}}\cr
\leq & \N {(\sum_j \pi (e_{jj}) \otimes \sum_i \xi_{ij}\xi_{ij}^*)^{1/2}}\cr
\leq & \sup_j \N {(\sum_i \xi_{ij}\xi_{ij}^*)^{1/2}}\cr
= & \N {\sum_{ij} e_{ii}\otimes e_{ij} \otimes \xi_{ij}}_{BR_n \otimes_{\min}
B(K_1)}}{.}$$
Taking $B(K_1) = M_n$ with $n \geq 1$ arbitrary, we obtain (recall (5.6))
$$\N {a_1}_{cb} \leq 1.$$
Similarly we have $\N {a_2}_{cb} \leq 1$.
This concludes the proof.

We now turn to a generalized version of Schur multipliers.

\n We consider an operator $T : M_n(B(H_1)) \to M_n(B(H))$ (where $H,H_1$ are
Hilbert spaces) of the special form 
$$T((x_{ij})) = (T_{ij} (x_{ij})) \leqno (5.13)$$
where $T_{ij}: B(H_1) \to B(H)$ are operators.

\n{\bf Remark 5.6}. Assume $\N {T}_{cb}\leq 1$. Then there is a Hilbert space
$K$, a representation $\pi : B(H_1) \to B(K)$ and operators $x_i,y_j : H \to
K$ such that for all $i,j$
$$\forall~x \in B(H_1)~~~T_{ij}(x) = x_i^* \pi(x) y_j~{\rm and} ~\N {x_i}
\leq 1, \N {y_j} \leq 1.\leqno (5.14)$$
Conversely, it is clear that (5.14) implies $\N {T}_{cb} \leq 1$.

\n This statement generalizes Proposition 1.1 to the present setting.

Such a statement is a simple consequence of the factorization theorem of c.b.
maps (\cf [Pa]) and of the particular form (5.13) of the map $T$.

\proclaim Corollary 5.7. Let $T : M_n(B(H_1)) \to M_n(B(H))$ be an operator of
the form (5.13) (generalized Schur multiplier). As above, let $G =
\tore^{n^2}$ and let $\mu$ be the normalized Haar measure on $G$. Let $T_z :
M_n (B(H_1)) \to M_n(B(H))$ be the operator defined by 
$$\forall~z = (z_{ij})\in G~~~~~T_z((x_{ij})) = (z_{ij}T_{ij} (x_{ij})).$$
Assume $\int \N {T_z}_{cb}d\mu(z) < 1$. 
Then there is a decomposition $T = \alpha_1 + \alpha_2$ where $\alpha_1$ and
$\alpha_2$ are each of the form (5.13) and moreover   there are
operators
$$\tilde \alpha_1 : BR_n \otimes_{\min}B(H_1) \to M_n(B(H))$$
and
$$\tilde \alpha_2 : BC_n \otimes_{\min}B(H_1) \to M_n(B(H))$$
satisfying
$$\N {\tilde \alpha_1}_{cb}\leq 1,\N {\tilde \alpha_2}_{cb} \leq 1,$$
and such that
$$\alpha_1 = \tilde \alpha_1(J_1 \otimes I_{B(H_1)})~~~~~\alpha_2 = \tilde
\alpha_2(J_2 \otimes I_{B(H_1)}).\leqno (5.15)$$
Conversely, if such a decomposition holds then necessarily
$$\int \N {T_z}_{cb}d\mu(z) \leq \sup_{z \in G}\N {T_z}_{cb}\leq 2.$$ 

\n{\bf Proof}. Let us define again $\tilde T_{ij} : M_n(B(H_1)) \to M_n(B(H))$ by
the identity
$$T_z = \sum_{ij}z_{ij} \tilde T_{ij}.$$
Then, by Theorem 5.1 there are a Hilbert space $\tilde K$ and a representation
$\tilde \pi : M_n(B(H_1)) \to B(\tilde K)$ together with operators
$V,W,V^{ij},W^{ij}$ from $\ell_2^n(H)$ into $\tilde K$ such that
$$\forall~\xi \in M_n(B(H_1))~~~~\tilde T_{ij}(\xi) = V^{ij*} \tilde \pi(\xi)
W + V^* \tilde \pi(\xi) W^{ij*}\leqno (5.16)$$
and such that
$$\N {V}\leq 1,\N {W}\leq 1, \N {\sum_{ij} V^{ij*}~V^{ij}}\leq 1,\N {\sum_{ij}
W^{ij*}~W^{ij}}\leq 1.\leqno (5.17)$$
By standard arguments, we can assume w.l.o.g. that $\tilde K = \ell_2^n(K)$ for
some Hilbert space $K$ and that $\tilde \pi : M_n(B(H_1)) \to B(\tilde K) =
M_n(B(K))$ is of the form $\tilde \pi = I_{M_n} \otimes \pi$ for some
representation $\pi : B(H_1) \to B(K)$. Once we identify $\tilde K$ with
$\ell_2^n(K)$ we may identify each of $V,W,V^{ij},W^{ij}$ with an $n \times n$
matrix of operators from $H$ into $K$.

Thus we identify $V$ with $(V(k,\ell))_{k,\ell \leq n}~~V^{ij}$ with $(V^{ij}
(k,\ell))_{k,\ell \leq n}$, and so on.

Let now $x$ be arbitrary in $B(H_1)$. We have by (5.16)
$$\eqalign{
T_{ij}(x)&= (\tilde T_{ij}(e_{ij} \otimes x))_{ij}\cr
& = [V^{ij*} (e_{ij}\otimes \pi(x)) W + V^*(e_{ij}\otimes \pi(x))
W^{ij}]_{ij}}$$
hence
$$T_{ij}(x) = V^{ij}(i,i)^* \pi(x) W(j,j) + V(i,i)^* \pi(x)
W^{ij}(j,j).\leqno(5.18)$$ By (5.17) we have 
$$\N {W(j,j)}\leq 1, \N {V(i,i)}\leq 1.$$
Moreover, we claim that (5.17) implies
$$\sup_i \N {\sum_j V^{ij}(i,i)^* V^{ij}(i,i)} \leq 1~{\rm and}~\sup_j \N
{\sum_i W^{ij} (j,j)^* W^{ij} (j,j)}\leq 1.\leqno(5.19)$$
Indeed, (5.17) implies that for each $\ell = 1,2,...,n$
$$\N {\sum_{ij} V^{ij} (\ell,\ell)^* V^{ij} (\ell,\ell)}\leq 1$$
hence a fortiori for each $i$ (taking $\ell=i)$
$$\N {\sum_j V^{ij} (i,i)^* V^{ij}(i,i)} \leq 1.$$
The other estimate is proved similarly, hence the above claim.

Finally, let
$$\alpha_1((x_{ij})) = (\alpha_{ij}^1 (x_{ij}))~{\rm and}~\alpha_2 ((x_{ij}))
= (\alpha_{ij}^2(x_{ij}))$$
where we define for all $x$ in $B(H_1)$
$$\alpha_{ij}^1(x) =
V(i,i)^* \pi(x) W^{ij}(j,j)~{\rm and}~\alpha_{ij}^2(x) = V^{ij}(i,i)^* \pi(x) W(j,j).$$
Clearly we can write (5.15) for some uniquely defined operators $\tilde
\alpha_1$ and $\tilde \alpha_2$ as in Corollary 5.7. We have 
$$\forall~x \in B(H_1)~~~\tilde \alpha_1(e_{jj} \otimes e_{ij} \otimes x) =
e_{ij} \otimes \alpha_{ij}^1(x)$$
$${\rm and}~~~~~~~~~~~~~~~\tilde \alpha_2 (e_{ii}\otimes e_{ij} \otimes x) =
e_{ij}\otimes \alpha_{ij}^2 (x).$$
Finally, it remains to check that $\tilde \alpha_1$ and $\tilde \alpha_2$
acting on the spaces indicated in Corollary 5.7 are of c.b. norm at most 1.

Let $y_1 \in BR_n \otimes_{\min} B(H_1)$ be of norm $\leq 1$. Let $y_1 = \sum
\limits_i e_{ii}\otimes \sum \limits_j e_{ij}\otimes x_{ij}$ with $x_{ij} \in
B(H_1)$. We have
$$\tilde \alpha_1(y_1) = \sum_{ij} e_{ij}\otimes V(i,i)^* \pi(x_{ij})
W^{ij}(j,j).$$
Observe that for all $y_{ij}$ in $B(H)$
$$\N {\sum e_{ij} \otimes y_{ij}}_{M_n(B(H))} = \N {\sum e_{ij} \otimes e_{ij}
\otimes y_{ij}}_{M_n(M_n(B(H))}$$
Now let $y_{ij}= V(i,i)^* \pi(x_{ij}) W^{ij} (j,j).$We have
$$\sum_{ij} e_{ij}\otimes e_{ij} \otimes y_{ij} = (\sum_{ij} e_{ii}\otimes
e_{ij}\otimes V(i,i)^* \pi(x_{ij})).(\sum_{ij} e_{ij}\otimes e_{jj}\otimes
W^{ij}(j,j)).$$
Hence we have
$$\N {\tilde \alpha_1(y_1)} \leq \N {\sum_{ij} e_{ii} \otimes e_{ij} \otimes
V(i,i)^* \pi(x_{ij})}.\N {\sum_{ij} e_{ij} \otimes e_{jj} \otimes W^{ij}(j,j)}$$
hence by (5.6)', (5.6)'' and (5.19).
$$ {
\N {\tilde \alpha_1(y_1)}  \leq \sup_i (\N {V(i,i)}\N {(\sum_j x_{ij}
x_{ij}^*)^{1/2}})  \leq \N {y_1}_{BR_n \otimes_{\min} B(H)}}{ }$$ 
This shows that $\N {\tilde \alpha_1}_{cb} \leq 1$. Similarly we have $\N
{\tilde \alpha_2}_{cb}\leq 1$. This concludes the proof.

\n{\bf Final Remark:} Recently, C.Le Merdy has shown that Theorem 5.1 and
corollary 5.7 remain valid if $B(H_1)$ is replaced by an arbitrary
$C^*$-algebra $A\subset B(H_1)$.
Indeed, he has proved (cf. [LeM]) that any bounded analytic function with values
in the space $CB(A,B(H))$ (i.e. the space of all   c.b. maps from $A$ into
$B(H)$) can be extended to a bounded analytic function (with the same
$H^\infty$ norm) with values in the space $CB(B(H_1),B(H))$. In other words,
there is a way   to extend c. b. maps
from $A$ into
$B(H)$ to c. b. maps defined on the whole of $B(H_1)$ which preserves
analyticity.
This extends a result
due to Haagerup and the author corresponding to the particular case when $H$ is
of dimension 1. Using Le Merdy's result and Remark 4.4 it is rather easy to adapt
the proof of Theorem 5.1 (or corollary 5.7) with  $B(H_1)$  replaced by any  
$C^*$-subalgebra $A\subset B(H_1)$.

\vfill\eject
\centerline {\bf REFERENCES}
\bg 
\bg 
\bg 
\ref AO & Akemann C. and Ostrand P.& Computing norms in group $C^*$-algebras. &
Amer. J. Math.& 98 (1976), 1015-1047.

\ref Be & Bennett G.&  Schur multipliers. & Duke Math.J.&
44(1977) 603-639.

\ref BP & Blecher D. and Paulsen V. & Tensor products of operator spaces  & J. Fuct. Anal. &
99 (1991) 262-292.

\ref B1 & Bo$\dot{z}$ejko M. & Remarks on Herz-Schur multipliers on free groups.
& Mat. Ann. & 258 (1981), 11-15.

\ref B2 & Bo$\dot{z}$ejko M. & Littlewood functions, Hankel multipliers and power
bounded operators on a Hilbert space. & Colloquium Math. & 51 (1987), 35-42.

\ref B3 & Bo$\dot{z}$ejko M. & Positive definite bounded matrices and a
characterization of amenable groups. & Proc. A.M.S. & 95 (1985), 357-360.

\ref BF & Bo$\dot{z}$ejko M. and Fendler G. & Herz-Schur multipliers and
completely bounded multipliers of the Fourier algebra of a locally compact
group. & Boll. Unione Mat. Ital. & (6) 3-A (1984), 297-302.

\ref DCH & de Canni\`ere J. and Haagerup U. & Multipliers of the Fourier algebras
of some simple Lie groups and their discrete subgroups.& Amer. J. Math. & 107
(1985), 455-500.

\ref C & Chevet S. & S\'eries de variables al\'eatoires gaussiennes  valeurs dans
$E \otimes_\eps F$. Application aux espaces de Wiener abstraits. Expos\'e XIX. &
S\'eminaire sur la g\'eom\'etrie des espaces de Banach. Ecole Polytechnique,
Palaiseau, & 1977-1978.

\ref D & Dean D. & The equation $L(E,X^{**})=L(E,X)^{**}$ and the principle of local reflexivity . &
Proc.A.M.S. & 40 (1973) 146-148.

\ref FTP & Figa-Talamanca A. and Picardello M. && Harmonic Analysis on Free
groups. & Marcel Dekker, New-York, 1983.

\ref GL & Gilbert J. and Leih T.& Factorization, Tensor
Products and Bilinear forms in Banach space Theory, in
Notes in Banach spaces,  (edited by E. Lacey). & University
of Texas Press, Austin &1980.

\ref G & Grothendieck A.& R\'esum\'e de la th\'eorie m\'etrique des produits tensoriels
topologiques. & Boll.. Soc. Mat. S$\tilde{a}$o-Paulo & 8 (1956), 1-79.

\ref H1 & Haagerup U. & An example of a non-nuclear
$C^*$-algebra which has the metric approximation property.
& Inventiones Mat. & 50 (1979), 279-293.

\ref H2 & Haagerup U. & Injectivity and decomposition of
completely bounded maps. in "Operator algebras and their
connection with topology and ergodic theory", & Springer
Lecture Notes in Math. & 1132 (1985) 91-116.

\ref H3 & Haagerup U. & Decomposition of completely
bounded maps on operator algebras.& Unpublished
manuscript.& Sept.1980.

\ref HP & Haagerup U. and Pisier G. &  Bounded linear operators
between $C^*$-algebras. & In preparation .&

\ref Ka & Kahane J.P. && Some randon
series of functions. Heath Math. Monograph. & 1968. New
edition Cambridge Univ. Press, 1985.

\ref K & Krivine J.L. & Th\'eor\`emes de factorisation
dans les espaces r\'eticul\'es.  & S\'eminaire
Maurey-Schwartz 73-74 Ecole Polytechnique, Paris &  
Expos\'es 22-23.

\ref Kw & Kwapie\'n S. & On operators factorizable
through $L_p$-space. & Bull. Soc. Math. France M\'emoire &
31-32 (1972) 215-225.

\ref L1 & Leinert M. & Faltungsoperatoren auf gewissen diskreten Gruppen. &
Studia Math. & 52, (1974), 149-158.

\ref L2 & Leinert M. & Abschtzung von Normen gewisser Matrizen und eine
Anwendung. & Math. Ann. & 240, (1979), 13-19.

\ref LeM & Le Merdy  C. & Analytic factorizatins and completely bounded maps. &
Preprint & To appear.

\ref LP & Lindalh L. and Poulsen F. && Thin sets in
Harmonic Analysis. & Marcel Dekker.New-York. 1971.

\ref LT & Lindenstrauss J. and Tzafriri L.&& Classical
Banach spaces II. & Springer Verlag Berlin, 1979.

\ref LR &  Lopez J. and Ross K. && Sidon sets. & Marcel
Dekker.New-York. 1975.

\ref LPP & Lust-Piquard F. and Pisier G. & Non commutative Khintchine and Paley
inequalities. & Arkiv fr Mat. & 29 (1991), 241-260.

\ref Pa & Paulsen V. && Completely bounded maps and dilations. Pitman Research.&
Notes in Math. 146, Longman, Wiley, New-York, 1986.

\ref P1 & Pisier G. && Factorization of linear operators and the Geometry of
Banach spaces.& CBMS (Regional conferences of the A.M.S.) n\up o 60, (1986),
Reprinted with corrections 1987.

\ref P2 & Pisier G. & Factorization of operator valued
analytic functions.& Advances in Math. & 93 (1992) 61-125.

\ref P3 & Pisier G. & Random series of trace class
operators. & Proceedings IV\up{o} C.L.A.P.E.M. (Mexico
Sept. 90). & To appear.

\ref P4 & Pisier G. & Completely bounded maps between sets of Banach space operators. &
Indiana Univ.
Math. J. & 39 (1990) 251-277. 

\ref R & Rudin W. && Fourier analysis on groups.&
Interscience. New York, 1962.

\ref V1 & Varopoulos N. & On an inequality of von Neumann and an application of
the metric theory of tensor products to Operators Theory. & J. Funct. Anal. &
16 (1974), 83-100.

\ref V2 & Varopoulos N. & Tensor algebras over discrete spaces. & J. Funct.
Anal. & 3 (1969), 321-335.

\ref W & Wysoczanski J. & Characterization of amenable groups and the
Littlewood functions on free groups.& Colloquium Math. & 55 (1988), 261-265.

\vskip12pt

Texas A. and M. University

College Station, TX 77843, U. S. A.

and

Universit\'e Paris 6

Equipe d'Analyse, Bo\^\i te 186,
 
75230 Paris Cedex 05, France

\end

\def\qed{{\vrule height7pt width7pt
depth0pt}\par\bigskip}
\def\p{{\bf Proof:\ }}

\def\p.{{\bf Proof:\ }}
(ii) $\Rightarrow$ (iii). We start by recalling that for any finitely supported
function $f : G \to B(H)$ we have the elementary inequality
$$\max \left\{\N {\left(\sum f(x)^* f(x)\right)^{1/2}},\N {\left(\sum f(x)
f(x)^*\right)^{1/2}}\right\} \leq \N {\sum \lambda(x) \otimes
f(x)}_{B(\ell_2(G,H))}.\leqno(3.3)$$
Now assume (ii). We have then by (3.3) if $\sup \limits_x |\eps(x)|\leq 1$
$$\N {\sum_{x \in \Lambda} \eps(x) \lambda(x) \otimes a(x)} \leq C \N {\sum_{x
\in G} \lambda(x) \otimes f(x)},$$
hence the multiplier of $C_\lambda^*(G)$ defined by $\eps$ is completely
bounded with norm $\leq C$. This proves (ii) $\Rightarrow$ (iii). 

(iii) $\Leftrightarrow$ (iv) follows from Proposition 1.2 above.

Finally   (iii) $\Rightarrow$ (ii) follows from, we have by duality for any
$\eps$ supported by $\Lambda$ with $\sup \limits_{x \in \Lambda} |\eps(x)|\leq
1$ and for any $f$ in $A(G,H \hat \otimes H)$.
$$\N {(\eps.f)}_{A(G,H \hat \otimes H)}\leq C \N {f}_{A(G,H \hat \otimes H)}$$
In particular, this holds for all the choices of signs $\eps(x)$ with
$|\eps(x)| = 1$, so that by the Theorem 0.1. in [LPP] (see Remark 3.2 above)
taking the average over all the choices of signs (or equivalently assume
$\Lambda$ countable, say $\Lambda = \{x_k|k \in \nat\}$, replace
$(\eps(x_k))_k$ by $(e^{i 3^k t})_k$ and average over $t$), we conclude that
$\Lambda$ is an $L$-set. This proves (iv) $\Rightarrow$ (i) and concludes the
proof.

 Indeed, by Proposition 3.3 the function $\ph$ which is the
indication function of $\Lambda$ satisfies property (i) hence (iii) and (iv) in
Theorem 0.1, hence properties (ii) and (iii) in Theorem 3.4 follow. Conversely,
we would give a simple argument for (iii) $\Rightarrow$ (i) but again this
follows from Proposition (3.3) and Theorem 0.1. as above.
Indeed, if $\ph$ satisfies (iii) in Theorem 0.1 then
Theorem 0.2 implies that there is a constant $C$ such that
for all finite subsets $S,T$ of $G$ and for all bounded
functions $\eps : G \to \tore$ with $|\eps (x)| \leq 1$ for
all $x$ in $G$, we have  $$\N {\sum_{s,t \in S \times  T}
\eps(st) \ph(st) e_s \otimes e_t}_{\ell_\infty (S) \hat
\otimes \ell_\infty(T)} \leq C$$ Clearly (recall Remark
1.4) this implies that the function $(s,t) \to \eps(st)
\ph(st)$ is a c.b. Schur multiplier of
$B(\ell_2(S),\ell_2(T))$ with norm $\leq C$. Since this
holds for all finite subsets $S,T$ (and since $\bigcup
\limits_S \ell_2(S)$ is dense in $\ell_2(G))$ the function
$(s,t) \to \eps(st) \ph(st)$ must also be a c.b. Schur
multiplier of $B(\ell_2(G),\ell_2(G))$ with norm $\leq C$.
Restricting now to the subspace $C_\lambda^*(G)$ (and using
the second assertion in Proposition 1.2) we find that
$\eps\ph$ is a c.b. multiplier, which proves (iii)
$\Rightarrow$ (i).

\end
\n{\bf Remark 2.5}. Let $G$ be a discrete group. Let us denote by $B_2(G)$ the
Banach space of all c.b. (= Herz Schur) multipliers on $G$. It is well known
(see [?]) that $B(G) = B_2(G)$ iff $G$ is amenable. It follows from Theorem
0.1 that if $B_2(G)$ is isomorphic (as a Banach space) to $B(G)$ or if
$B_2(G)$ is isomorphic to a subspace of $B(G)$ then $G$ is amenable. Indeed by
[TJ] $B(G)$ is of cotype 2 for any group, hence if $B_2(G)$ is isomorphic to a
subspace of $B(G)$ then $B_2(G)$ itself is of cotype 2, hence by Theorem 0.1
the space $T_2(G)$ (of all functions satisfying the third property in Theorem
0.1) is of cotype 2, hence $T_2(G) \subset \ell_2(G)$, and we conclude by [W]
that $G$ is amenable. This suggests the following question : If $B_2(G)$ is of
cotype $q$ for some $q < \infty$, is $G$ amenable ?

\n{\bf Final Remark:} Recently, C.Le Merdy has shown that Theorem 5.1 and
corollary 5.7 remain valid if $B(H_1)$ is replaced by an arbitrary
$C^*$-algebra $A\subset B(H_1)$.
Indeed, he has proved (cf. [LeM]) that any bounded analytic function with values
in the space $CB(A,B(H))$ (i.e. the space of all   c.b. maps from $A$ into
$B(H)$) can be extended to a bounded analytic function (with the same
$H^\infty$ norm) with values in the space $CB(B(H_1),B(H))$. In other words,
there is a way   to extend c. b. maps
from $A$ into
$B(H)$ to c. b. maps defined on the whole of $B(H_1)$ which preserves
analyticity.
Le Merdy's theorem extends the result
of [HP] corresponding to the particular case when $H$ is of dimension 1.
Using Le Merdy's result and Remark 4.4 it is rather easy to adapt the proof of
Theorem 5.1 (or corollary 5.7) with  $B(H_1)$  replaced by any  
$C^*$-subalgebra $A\subset B(H_1)$.
\bye